\documentclass[dvips]{amsart}
\usepackage{amssymb,amsxtra,latexsym,pictex,epsfig,psfrag}

\newcommand{\vanish}[1]{}
\newcommand{\F}{\mathcal{F}}
\newcommand{\G}{\mathcal{G}}
\newcommand{\A}{\mathcal{A}}

\newcommand{\C}{\mathcal{C}}
\newcommand{\T}{\mathcal{S}}
\newcommand{\D}{\Delta}

\newcommand{\Z}{\mathbb{Z}}
\newcommand{\fq}{\mathbb{F}_q}
\newcommand{\abs}[1]{\lvert#1\rvert}

\newcommand{\zerohat}{\hat{0}}
\newcommand{\onehat}{\hat{1}}
\newcommand{\supp}{\operatorname{supp}}

\newcommand{\Des}{\operatorname{Des}}

\newcommand{\Stab}{\operatorname{Stab}}
\newcommand{\onto}{\twoheadrightarrow}
\newcommand{\into}{\hookrightarrow}

\theoremstyle{plain}
\newtheorem*{corollary}{Corollary}
\newtheorem{lemma}{Lemma}
\newtheorem{theorem}{Theorem}

\newtheorem{proposition}{Proposition}

\theoremstyle{definition}

\theoremstyle{remark}
\newtheorem*{remark}{Remark}
\newtheorem{example}{Example}

\newcommand{\abra}{abramenko94:_walls_coxet}
\newcommand{\abels}{abels91}
\newcommand{\bjorner}{bjoerner84:_some_coxet_tits}
\newcommand{\bjor}{bjoerner80:_shell_cohen_macaul}
\newcommand{\bw}{bjoerner97:_shell}
\newcommand{\riffle}{bayerdiaconis92:_trail}
\newcommand{\bid}{bidigare97:_hyper}

\newcommand{\bhr}{bidigarehanlon98}
\newcommand{\bbd}{billera99:_random}

\newcommand{\bir}{birkhoff79:_lattic}
\newcommand{\red}{bjoerner93:_om}
\newcommand{\brady}{brady01}
\newcommand{\brown}{brown89:_build}
\newcommand{\ken}{brown00:_semig_markov}
\newcommand{\bd}{browndiaconis:_random}

\newcommand{\grillet}{grillet95:_semig}
\newcommand{\gruna}{gruenbaum72:_arran}

\newcommand{\muhlherr}{muehlherr90}

\newcommand{\klein}{klein-barmen40:_uber_veral_verban}
\newcommand{\me}{math.PR/0108094}
\newcommand{\ot}{orlikterao92:_book}
\newcommand{\petrich}{petrich71}
\newcommand{\petrichbook}{petrich77:_lectur}

\newcommand{\ronan}{ronan89:_lectur}
\newcommand{\dress}{dress87:_gated}
\newcommand{\scharlau}{scharlau85:_metric}
\newcommand{\schar}{scharlau95:_build}
\newcommand{\schmid}{muehlherr95:_stein}
\newcommand{\schutz}{schuetzenberger47:_sur}

\newcommand{\solomon}{solomon76:_mackey}

\newcommand{\stanley}{stanley97:_enumer1}
\newcommand{\stanleysuper}{stanley72:_super}
\newcommand{\tits}{tits74:_build_bn}

\newcommand{\zie}{ziegler95:_lectures}

\begin{document}
\title{Projections, shellings and duality}
\author{Swapneel Mahajan}

\address{Department of Mathematics\\
Cornell University\\
Ithaca, NY 14853}
\email{swapneel@math.cornell.edu}


\begin{abstract}
Projection maps which appear in the theory of buildings
and oriented matroids are closely related to the notion of shellability.
This was first observed by Bj{\"o}rner~\cite{\bjorner}.
In this paper, we give an axiomatic treatment of either concept and
show their equivalence. We also axiomatize duality in this setting.
As applications of these ideas,
we prove a duality theorem on buildings and give a geometric
interpretation of the flag $h$ vector. 
The former may be regarded as a $q$-analogue of the Dehn-Sommerville equations.
We also briefly discuss the connection with the random walks
introduced by Bidigare, Hanlon and Rockmore~\cite{\bhr}.
\end{abstract}

\maketitle

\section{Introduction}

Projection maps appeared independently in the theory of buildings
developed by Tits~\cite{\tits} and in the theory of oriented
matroids~\cite{\red}. 
Given the importance of these maps in either theory,
it is reasonable to try to formulate them axiomatically.
The close connection between these maps and the notion of shellability
was first observed by Bj{\"o}rner~\cite{\bjorner}.
In this paper, we study projection maps axiomatically keeping the
viewpoint of~\cite{\bjorner}. This leads us to the notion of many
``compatible'' shellings rather than a single shelling.
To complete the picture, we also consider restriction maps that are
useful to keep track of a shelling. 
For ideas closely related to this paper, see~\cite{\scharlau,\dress,\muhlherr}.
The main result of this paper (Theorem~\ref{t:axiom}) may be
informally stated as 
$$\text{Projection maps} \Longleftrightarrow \text{Compatible shellings}.$$

Thus far, people seem to be interested in the question 
``Is a complex shellable?''
The next question that should be asked is
``How do we study the space of all shellings of a complex?''
As explained in the previous paragraph,
projection maps give us an approach to this hard problem.
We mention that though compatible shellings form a fairly large class
of shellings, they are far from giving all possible shellings of a
complex.

For a labelled shellable complex $\D$, Bj{\"o}rner~\cite[Theorem
1.6]{\bjorner} gives a simple geometric 
interpretation of the flag $h$ vector.
In our case, $\D$ has many ``compatible'' shellings rather than a
single shelling. This gives us a dual
interpretation of the flag $h$ vector which is easier to visualise.
It is particularly nice when the projection maps
satisfy some commutativity relations. 
These relations are related to the uniformity of the stationary distribution of
certain random walks introduced by Bidigare, Hanlon and Rockmore~\cite{\bhr}.
In this context, we generalise a result on uniform stationary
distributions obtained in~\cite{\bbd}.

Next we formalise the notion of duality by adding an ``opposite''
axiom to our axiomatic setup. It is strong enough to imply the
Dehn-Sommerville equations. 
And as one expects, this axiom can hold only for a simplicial complex homotopy
equivalent to a sphere. 
From the viewpoint of shellings, the
relevant concept is that of shelling reversal, i.e., when is the
reverse of a shelling again a shelling?

The Solomon-Tits theorem says that a (thick) spherical building $\D$ has
the homotopy type of a (non-trivial) wedge of spheres.
Hence by the observation in the previous paragraph, 
$\D$ cannot satisfy the opposite axiom. 
However $\D$ has a remarkable duality which can be expressed in terms of
the flag $h$ vector (Theorem~\ref{t:duality}).
This may be regarded as a $q$-analogue of the Dehn-Sommerville equations.
 
\subsection*{Organization of the paper}
The next two sections provide basic definitions and motivating
examples.
The axiomatic theory involving projections maps, shellings and
restriction maps is presented in Section~\ref{s:axiom} and later in a
more intuitive metric setup in Section~\ref{s:metric}.
In Section~\ref{s:tss} we explain the connection with 
the flag $h$ vector.
The connection with random walks is explained in Section~\ref{s:comm}.
In the next two sections, we study the notion of duality
first for thin complexes and then for buildings.
In the final section, we outline some problems for further study.

\section{Projection maps} \label{s:proj}

In this section, we give four examples of the mini-theory that we will
present in Section~\ref{s:axiom}.
We will return to these examples again in Section~\ref{s:metric},
where we will show that they are indeed examples of our theory.
Terms like gallery connected, convex, etc. that we freely use here
are also explained in that section.
For now, we present the examples from the viewpoint of projection maps. 
They fit into the framework of LRBs (non-associative in
general) as indicated in the table.
\\

\begin{center}
\begin{tabular}{|c|c|} \hline 
LRB & Non-associative LRB \\ \hline \hline
Hyperplane arrangements & ??? \\ \hline
Reflection arrangements & Buildings \\ \hline
Distributive lattices & Modular lattices \\ \hline 
\end{tabular}
\end{center}

\vspace{0.3 cm}
\noindent
The question marks say that a building-like analogue for the more
general case of hyperplane arrangements is unknown.  
Though LRBs in general are not examples of our theory they come quite
close as we will see in Section~\ref{s:LRB}; and hence can potentially
give more examples.
Keeping this in mind, we start with a brief review of LRBs.

\subsection{Left regular bands}

Let $S$ be a semigroup (finite, with identity).  
A \emph{left-regular band}, or LRB, is a semigroup $S$ that satisfies
the identities
\begin{equation} \tag{D}
x^2=x \quad \text{and} \quad xyx=xy
\end{equation}
for all $x,y\in S$.  We call (D) the ``deletion property'', because it
admits the following restatement:  Whenever we have a product $x_1
x_2\cdots x_n$ in~$S$, we can delete any factor that has occurred
earlier without changing the value of the product.  

Alternatively, one can say that $S$ is a \emph{LRB} if there are a
lattice $L$ and a  
surjection $\supp\colon S\onto L$ satisfying
\begin{equation} \label{e:homom}
\supp xy = \supp x \vee \supp y
\end{equation}
and
\begin{equation} \label{e:delete}
xy=x \quad\text{if } \supp y\le \supp x.
\end{equation}
Here $\vee$ denotes the join operation (least upper bound) in $L$.  
A good reference for LRBs is \cite{\ken}. It explains the equivalence
of the above two definitions and also contains plenty of examples.
More information about
LRBs can be found in \cite{\grillet,\petrich,\petrichbook}.  Early
references to the identity $xyx=xy$ are \cite{\klein,\schutz}.

We can also define a partial order on a LRB by
setting
\begin{equation} \label{e:poset1}
x\le y \iff xy=y.
\end{equation}

\vanish{
Note that left multiplication by $x$ is a \emph{projection}
(idempotent operator) mapping $S$ onto $S_{\ge x} = \{y\in S : y\ge
x\}$.  The latter is a LRB in its own right, the
associated lattice being the interval $[X,\onehat]$ in~$L$, where
$X=\supp x$ and $\onehat$ is the largest element of~$L$.  Note also
that $S_{\ge x}$ depends only on $X$, up to isomorphism.  
We may therefore write $S_{\ge X}$ instead
of $S_{\ge x}$.
We will need the subLRB $S_{\ge X}$ briefly in Section~\ref{s:comm}.
}

\noindent
We now switch to a slightly different notation which we will be using
for the most part.
We denote a LRB by $\F$ and call elements of $\F$ as faces.
A face $C$ is a \emph{chamber} if $CF=C$ for all $F \in \F$,
or equivalently, if it is maximal in the partial order on $\F$ 
specified by equation~\eqref{e:poset1}.
The set of chambers, which we denote $\C$, is an ideal in $\F$.
This gives us a map $\F \times \C \rightarrow \C$ that maps $(F,C)$ to $FC$.
We call $FC$ the \emph{projection} of $C$ on $F$.
In this paper, we are mainly interested in these projection maps
rather than the full product in $\F$.

\begin{example} \label{e:hfs}
Hyperplane arrangements:
The motivating example of a LRB is the poset of regions of a central
hyperplane arrangement (or more generally of an oriented matroid).
A good reference for this example is \cite[Appendix A]{\ken}.
More details can be found in
\cite{\bhr,\bbd,\red,\brown,\bd,\ot,\zie}.  
Briefly, a finite set of linear hyperplanes 
(i.e. hyperplanes passing through the origin)
in a real vector space $V$
divides $V$ into regions called \emph{chambers}.  These are polyhedral
sets, which have \emph{faces}.  The totality $\F$ of all the faces is a poset
under the face relation.  Less obviously, $\F$ admits a product,
making it a LRB.  
The lattice $L$ is the intersection lattice of the arrangement (or
more generally the underlying matroid).

The product can be described combinatorially by
encoding faces using sign sequences.
In fact one way to axiomatize an oriented matroid is in terms of sign
sequences and this product. 
The deletion property $(D)$ for $\F$ can be checked directly.
The product in $\F$ can also be described geometrically. For instance,
the projection $FC$ is the chamber closest to $C$ 
having $F$ as a face.

The case that is directly relevant to our theory is that of \emph{simplicial}
hyperplane arrangements, 
i.e., the chambers are simplicial cones and hence 
$\F$ is a simplicial complex.
As an interesting example, we mention
Coxeter complexes which arise from reflection arrangements.

Now let $\mathcal{D} \subseteq \C$ be a convex set of chambers and let
$\G$ be the set of faces of all the chambers in $\mathcal{D}$.
Then using the geometric description of the product in $\F$, for
instance, one can check that $\G$ is
also a LRB. We now give a concrete example of this type.
\end{example}

\begin{example}
Distributive lattices:
A good reference for this example is \cite[Section 4]{\ken}.
We will soon see that it generalises to the
case of modular lattices (Example~\ref{e:modular}). 
We start with the basic example of the Boolean lattice $\mathcal{B}_{n}$
 of rank $n$ 
consisting of all subsets of an $n$-set ordered under inclusion.
Its flag (order) complex $\D(\mathcal{B}_{n})$ 
is the Coxeter complex of type $A_{n-1}$ and
corresponds to the braid arrangement. 
Hence the set of faces of $\D(\mathcal{B}_{n})$ is a LRB.
Moreover the product can be described entirely using the lattice
structure (meets and joins) of $\mathcal{B}_{n}$.

More generally, the set of faces $\F$ 
of the flag (order) complex $\D(M)$ of any distributive lattice
$M$ is also a LRB. 
A simple way to verify this is
to appeal to the well-known fact that $M$ can be embedded as a
sublattice of the Boolean lattice.  More geometrically,
Abels~\cite[Proposition~2.5]{\abels} has described a way of
constructing an embedding which makes the set of chambers in~$\D(M)$
(i.e., the maximal chains in~$M$) correspond to a convex set of
chambers in $\D(\mathcal{B}_{n})$, 
the complex of the braid arrangement.

\end{example}

\subsection{Non-associative LRBs} 

As one may expect, these
are the same as LRBs except that we no longer require
associativity. However, we do require that $xyx$ be well-defined for
all $x,y$; that is, $xyx=(xy)x=x(yx)$.
With this restriction, the deletion property makes sense.
In the non-associative setting, the second definition of LRBs
involving a lattice $L$ does not make much sense and we omit it.
Also transitivity of the relation $\leq$ defined by
equation~\eqref{e:poset1} is no longer automatic. So we have two
choices; either to take transitive closure of the relation or simply to impose
transitivity as an additional condition.
In the two examples that we consider, transitivity of the relation
$\leq$ is in fact automatic and hence we do not
pursue this issue further.

To avoid confusion later, we mention that the term LRB always means an
associative LRB. Whenever we want to include the non-associative case,
we will say so explicitly.

\begin{example}
Buildings:
A good reference for this example is \cite{\tits}; also
see~\cite{\brown,\schar}. 
Some of the terminology used here is explained at the beginning of
Section~\ref{s:metric}. 
Also the building of type $A_{n-1}$ is described briefly in
Example~\ref{e:dbuild} in Section~\ref{s:dbuild}.
Let $W$ be a Coxeter group and $\Sigma(W)$ its Coxeter complex.
Roughly
a building $\D$ of type $W$ is a union of subcomplexes $\Sigma$ (called
\emph{apartments}) which fit together nicely.
Each apartment $\Sigma$ is isomorphic to $\Sigma(W)$.
For any two simplices in $\D$, there is an apartment $\Sigma$
containing both of them.
As a simplicial complex, $\D$ is pure, labelled and gallery connected. 
Also each apartment is convex.
There is a $W$-valued distance function $\delta:\C \times \C
\rightarrow W$ that generalises the gallery distance.
Furthermore for any apartment $\Sigma$ and chamber $C \in \C$, 
there is a retraction 
$\rho=\rho_{\Sigma,C} : \D \rightarrow \Sigma$ satisfying
$\delta(C,\rho(D))=\delta(C,D)$ for any $D \in \C$.

We denote the set of faces by $\F$.
For $F,G \in \F$, we choose an apartment $\Sigma$ containing $F$ and $G$
and define $FG$ to be their product in $\Sigma$.
Since $\Sigma$ is a Coxeter complex, we know how to do this
(Example~\ref{e:hfs}).  
Furthermore it can be shown that the product does not depend on the
choice of $\Sigma$.
Hence this defines a product on the set of faces $\F$ of $\D$.
And it is compatible with the retraction $\rho=\rho_{\Sigma,C}$, namely,
$\rho(F) C = \rho(F C)$ for any face $F$ of $\D$.
\end{example}

\begin{example} \label{e:modular}
Modular lattices:
A good reference for this example is \cite{\abels} where you can find
proofs of all the facts that we state here.
The flag (order) complex $\D(M)$ of a modular lattice $M$ is a labelled
simplicial complex.
A face of $\D$ is a chain in $M$. 
As usual, we denote the set of faces by $\F$.
The chambers $\C$ of $\D$ are the maximal chains in $M$.
We define the map 
$\F \times \C \rightarrow \C$ as follows.
For $\F \in \F$ and $C \in \C$, define $FC$ to be the unique chamber
containing $F$ that is contained in the sublattice generated by $F$ and
$C$.
More generally, for $F,G \in \F$, the face $FG$ is the chain in $M$ obtained by
refining the chain $F$ by the chain $G$ (using meets and joins as in a
Jordan--H\"{o}lder product).

An interesting subclass of modular lattices is that of distributive
lattices that we discussed earlier. 
Another interesting fact is that in a modular lattice, 
the sublattice generated by any two chains is distributive; see
\cite[pg 66]{\bir}.
Hence modular lattices may be regarded as generalisations of buildings of
type $A_{n-1}$, with these distributive lattices playing the role of
apartments (and their convex subsets).

It may be possible to generalise this example in various
directions. For example, one may consider more general lattices like
supersolvable lattices~\cite{\stanleysuper}. 
It is also a challenge to find analogues of modular lattices that
generalise buildings of types other than $A_{n-1}$. 

\end{example}

\begin{remark}
For buildings, while there is always an apartment containing two given
chambers, there may not be an apartment containing three given
chambers.
Similarly for modular lattices, the sublattice generated by three
chains may not be distributive.
This is the basic reason why the product in buildings and modular
lattices is not associative. We do not know of any algebraic tools
relevant to the study of non-associative LRBs.
\end{remark}

\begin{remark}
All four examples that we discussed share some common geometric properties.
The LRB $\F$ (non-associative in the last two examples) 
is a simplicial (or polyhedral) complex.
(Later in the paper, we will call such a $\F$, a \emph{simplicial LRB}.)
Furthermore $\F$ is pure and gallery connected.
Geometrically the projection $FC$ is the chamber closest to $C$
having $F$ as a face.
Also $FG = \cap_{G<D,D \in \C} FD$ for $F,G \in \F$.
\end{remark}

\section{Shellings and restriction maps} \label{s:shelling}

Let $\D$ be a finite pure $d$-dimensional simplicial complex.
The term \emph{pure} means that all maximal simplices have the same
dimension. We will call the maximal simplices as chambers.
Let $\F$ be the set of all faces of $\D$ and let $\C$ be the set of
chambers of $\D$.
Also let $\D_{\geq F}$ be the simplicial subcomplex of $\D$ consisting
of all faces of the faces that contain $F$.
Similarly $\C_{\geq F}$ stands for the set of chambers of $\D_{\geq F}$.
Let the partial order $\leq$ denote the face relation; that is,
$F \leq G$ if $F$ is a face of $G$. 

We say that
$\D$ is \emph{shellable} if there is a linear order $\leq_S$ on the set of
chambers $\C$ of $\D$ such that 
for every $D \in \C$ except the first in the linear order, 
we have $D \cap (\cup_{E<_S D} E)$ is pure $(d-1)$-dimensional;
that is, the intersection of $D$ with the chambers that came before it
in the linear order $\leq_S$ is a non-empty union of certain facets of $D$.
Less formally, a shelling gives a systematic way to build $\D$ by
adding one chamber at a time. And the subcomplex obtained at 
each stage (and in particular the entire complex $\D$) is gallery connected.

\begin{figure}[hbt]
\centering
\begin{tabular}{c@{\qquad}c}
\mbox{\epsfig{file=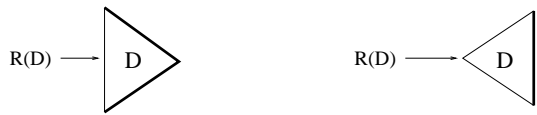,width = 7cm}}
\end{tabular}
\caption{The restriction map at work.}
\label{residue}
\end{figure}

To a shelling $\leq_S$, we can associate a \emph{restriction map}
$R : \C \rightarrow \F$ with
$R(D)$ defined to be the face of $D$ spanned by those vertices $v$ of $D$ for which
$D \setminus v \leq D \cap (\cup_{E<_S D} E)$.
Here $D \setminus v$ is the facet of $D$ that does not contain the
vertex $v$.
Figure~\ref{residue} shows two cases that can arise in the rank 3
case. The facets of $D$ that are shown by dark lines indicate the
intersection of $D$ with the chambers that came before it
in the shelling order.
Observe that $R(D) \leq F \leq D \Leftrightarrow F \ \text{shows up
for the first time when we adjoin} \ D$.
This shows that $\F$ can be expressed as a disjoint union 

\begin{equation} \label{e:disjoint}
\F = \sqcup_{D \in \C} \{ F \in \F \mid R(D) \leq F \leq D \}.
\end{equation}
Compare this statement with the restriction axiom $(R2)$ in
Section~\ref{subs:axioms}.
Also note that 
$R(D)=\emptyset$, the empty face of $\D$, if and only if 
$D$ is first in the linear order $\leq_S$.
To summarise, the restriction map $R$ gives us the local data at
every chamber associated with the shelling $\leq_S$. It does not have
enough information however to allow us to reconstruct $\leq_S$. 
This is because the essence of a shelling is really a partial order.
Here is what we can do. 
We can define a partial order $\leq_C$ on $\C$ as the transitive 
closure of the relation:
\begin{equation} \label{e:con}
E \leq_C D \ \text{if} \ R(E) \leq D.
\end{equation}
Observe that the chamber $C$ by which we indexed our partial order is
the one that occurred first in $\leq_S$. 
The motivation for this notation will be more clear when we study the
shelling axioms in Section~\ref{subs:axioms}. 
One can show directly that any linear extension of $\leq_C$ is a shelling of
$\D$ (compare this statement with the shelling axiom $(S2)$) and the
shelling $\leq_S$ that we started with is one of them.
We can also define a more refined partial order $\leq_C^r$ on $\C$ 
as the transitive closure of the relation:
\begin{equation} \label{e:rcon}
E \leq_C^r D \ \text{if} \ R(E) \leq D \ \text{and} \ E \ 
\text{is adjacent to} \ D.
\end{equation}
We will show (Lemma~\ref{l:comp}) that for a thin complex,
the partial orders $\leq_C$ and $\leq_C^r$ on $\C$ are identical.
The term \emph{thin} means that every facet is contained in exactly 2
chambers. 

Now suppose instead that we start not with a shelling order $\leq_S$
but only with a map $R : \C \rightarrow \F$
that satisfies equation~\eqref{e:disjoint}.  
Then can we say that $R$ is the restriction map of a shelling?
We know only a partial answer.
The first thing to do is to define a relation $\leq_C$
using equation~\eqref{e:con}.
(The chamber $C$ is the unique chamber satisfying $R(C) = \phi.)$
If the relation $\leq_C$ happens to be a partial order 
then it follows that any linear extension of
$\leq_C$ is a shelling of $\D$ 
and the restriction map of any of these is $R$.

The above statements will be justified 
in the course of proving Theorem~\ref{t:axiom}.
For future use, we record three useful results about shellable complexes.

\begin{lemma} \label{l:sh}
Let $F$ be any face of $\Delta$ and let $\leq_S$ be a shelling
of $\Delta$. Then this linear order when restricted to $\C_{\geq F}$
is a shelling of $\Delta_{\geq F}$.
In other words, the link $lk(F,\D)$ of any face $F$ in a shellable
complex $\D$ is again shellable with the induced order. 
\end{lemma}
\begin{proof}
Let $\leq_S$ be a shelling of $\D$ 
and let $\leq_S^{\prime}$ be its restriction to $\C_{\geq F}$.
By the definition of a shelling we have
\begin{equation} \tag{*}
D \cap (\cup_{E <_S D} E) = \cup_{G \in \F_D} G,
\end{equation}
where $\F_D$ is a subset of the set of facets of $D$.
Let $\F^{\prime}_D$ be the subset of $\F_D$ consisting of those facets that contain
$F$. The lemma follows from the following claim.

\medskip
\noindent
Claim: For $D \geq F$, we have $D \cap (\cup_{E <_S^{\prime} D} E) = \cup_{G \in \F^{\prime}_D}
G$.\\
Note that if $D$ is not the first element in $\leq_S^{\prime}$
then $F \subseteq LHS$.
In that case the set $\F^{\prime}_D$ will be non-empty as required.\\
Proof of the claim. $(\supseteq)$ Let $G \in \F^{\prime}_D$. 
By equation $(*)$, there is a $E \geq G$ such that $E <_S D$. 
Since $G \geq F$, it follows that $E \geq F$ and hence $E <_S^{\prime} D$.\\
$(\subseteq)$ Let $F^{\prime} \subseteq LHS$. We may assume that $F^{\prime} \geq F$.
Since otherwise we may replace $F^{\prime}$ by the face spanned by $F$
and $F^{\prime}$ which still belongs to the LHS.
Applying $(*)$, there is a $G \in \F_D$ such that $G \geq F^{\prime} \geq F$.
So $G \in \F^{\prime}_D$. Hence $F^{\prime} \subseteq RHS$.
\end{proof}

\begin{lemma} \label{l:comp}
Let $\D$ be a thin shellable complex 
with shelling $\leq_S$ and restriction map $R$.
Let $\leq_C$ and $\leq_C^r$ be the partial orders on $\C$ 
defined by equations~\eqref{e:con} and~\eqref{e:rcon} respectively.
Then the partial orders $\leq_C$ and $\leq_C^r$ on $\C$ are identical.
\end{lemma}
\begin{proof}

We only need to show that 
$E \leq_C D \Rightarrow E \leq^{r}_C D$,
or equivalently,  
$R(E) \leq D$ implies $E \leq^{r}_C D$.
We do this by constructing a gallery from $E$ to $D$ such that for
consecutive chambers $E^{\prime}$ and $E''$ in the gallery, we have
$R(E^{\prime}) \leq E''$.

\begin{figure}[htb]
\begin{center}
\psfrag{E}{\Huge $E$}
\psfrag{E1}{\Huge $E'$}
\psfrag{E2}{\Huge $E''$}
\psfrag{D}{\Huge $D$}
\psfrag{F}{\huge $F$}
\resizebox{!}{3cm}{\includegraphics{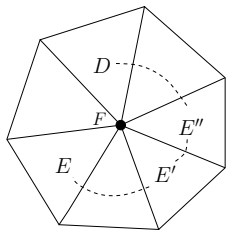}}
\end{center}
\caption{A gallery from $E$ to $D$ in $\D_{\geq F}$ where $F=R(E)$.}
\label{refine}
\end{figure}

\vanish{

\begin{figure}[hbt]
\centering
\begin{tabular}{c@{\qquad}c}
\mbox{\epsfig{file=refine.eps,height=3cm,width = 3cm}}
\end{tabular}
\caption{A gallery from $E$ to $D$ in $\D_{\geq F}$ where $F=R(E)$.}
\label{refine}
\end{figure}

}

\noindent
By Lemma~\ref{l:sh}, the shelling $\leq_S$ 
restricts to a shelling $\leq_S^{\prime}$ of $\D_{\geq R(E)}$. 
Hence $\D_{\geq R(E)}$ is gallery connected 
at every stage of the shelling $\leq_S^{\prime}$.
And $E$ is the starting chamber of the shelling. 
This allows us to choose a gallery from $E$ to $D$ in $\D_{\geq R(E)}$
consistent with the shelling $\leq_S^{\prime}$.
This automatically implies 
consistency with the original shelling $\leq_S$.
Let $E^{\prime}$ and $E''$ be two consecutive chambers
in the gallery 
such that $E^{\prime} <_S E''$.
Since $\D$ is thin, the common facet of 
$E^{\prime}$ and $E''$ appears for the first time in the shelling
$<_S$ when we adjoin $E^{\prime}$. 
Hence we get $R(E^{\prime}) \leq E''$.
This is exactly what we wanted to show.
Hence the partial orders $\leq_C$ and $\leq^{r}_C$
are identical.

\end{proof}

\begin{proposition} \label{p:shell}
\cite[Theorem 1.3]{\bjorner}
A finite $d$ dimensional shellable complex $\D$ is homotopy equivalent
to a wedge of $d$-spheres. The number of these spheres is given by
$\abs{\{D \mid R(D) = D\}}$, 
where $R$ is the restriction map associated with any shelling of $\D$.
\end{proposition}

Note that if $\abs{\{D \mid R(D) = D\}} = 1$
then the unique chamber $D$ for which $R(D) = D$ is the one that gets
shelled in the end.

\section{Axioms} \label{s:axiom}

In this section we provide an axiomatic setup that relates projection
maps, shellings and restriction maps. 
Motivation for some of the axioms was given in the previous two sections.
The axioms are somewhat abstract and you may want to simply glance
at them now.
They are split into three categories.
The first two axioms in each of the three categories really deal with a fixed
chamber $C \in \C$.
As we vary $C \in \C$, it is natural to impose some \emph{compatibility}
condition. This is the content of the third (compatibility) axiom.
It is easiest to swallow in the shelling case. We do not know of any
way to make it more palatable in the other two cases.

The main result of this section (Theorem~\ref{t:axiom})
says that for a pure simplicial
complex, the different sets of axioms are equivalent.
The way to pass from one set of axioms to another is explained
immediately after the statements of all the axioms.
The examples mentioned in Section~\ref{s:proj} are discussed
from this axiomatic viewpoint in Section~\ref{s:metric}.
It is a good idea to read that section before reading the proof of 
Theorem~\ref{t:axiom}.
Apart from serving as motivation, it will unable you to understand
the geometric content of every step in the proof.

Let $\D$ be a finite pure simplicial complex.
Let $\F$ be the set of all faces of $\D$ and let $\C$ be the set of
chambers of $\D$.
We always consider the empty face $\emptyset$ to be a face of $\D$.
Let the partial order $\leq$ denote the face relation; that is,
$F \leq G$ if $F$ is a face of $G$. 
Also let $\lessdot$ stand for the cover relation; that is,
$G \lessdot D$ if $G$ is a codimension $1$ face of $D$. 
If $D$ is a chamber and $G \lessdot D$ then $G$ is a facet of $\D$.

\subsection{The axioms} \label{subs:axioms}

Before we state the projection axioms, we need a definition.
Given a map $\F \times \C \rightarrow \C$, we say that 
$D$ and $E$ are \emph{weakly $C$ adjacent} if they have a common face
$F$ such that $FC=D$.
We may write this as 
$D \ \frac{F}{} \ E$.

Note that this is not a symmetric relation.
Also we do not require $F$ to be a facet.
In fact, $F$ could also be the empty face.
This explains the term ``weakly adjacent''.
A \emph{weak $C$ gallery} is
a sequence of chambers $C_1,\ldots,C_n$ 
such that consecutive chambers are weakly $C$ adjacent.
We write this as $C_1-\ldots-C_n$. 
If more precision is required then we write
$C_1\ \frac{F_1}{} \ C_2 \ \frac{F_2}{}\ldots\frac{F_{n-1}}{} \ C_n=D,$
where $F_1,\ldots,F_{n-1}$ 
is a sequence of faces 
such that, $F_i C = C_i$ and $F_i \leq
C_i,C_{i+1}$. 
\begin{figure}[hbt]
\centering
\begin{tabular}{c@{\qquad}c}
\mbox{\epsfig{file=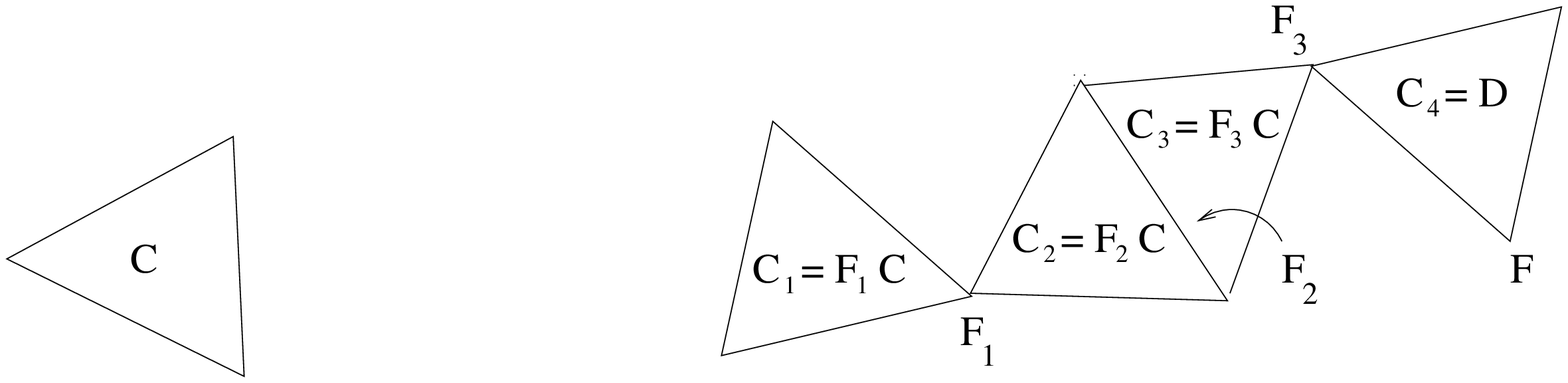,height=2.7cm,width = 10cm}}
\end{tabular}
\caption{A weak $C$ gallery $C_1 - C_2 - C_3 - C_4$.}
\label{gallery}
\end{figure}
\noindent
It will be seen from axiom $(P1)(ii)$ that $\emptyset C = C$.
Hence $C$ can always be tagged on as the first element of any weak $C$
gallery. 

\subsection*{Projection axioms $(P)$}
For every $F \in \F$, there is a projection map $\C \rightarrow \C$ which we
write $C \mapsto FC$ that satisfies

\begin{itemize}
\item[(P1)]

(i) If $FC=D$ then $F \leq D$.\\
(ii) If $F \leq C$ then $FC=C$.\\
(iii) If $FC=D$ and $F \leq G \leq D$ then $GC=D$.

\item[(P2)]

If $F \leq D$ and $GC=D$ for all $F \leq G \lessdot D$ then $FC=D$.
As a special case when $F=\emptyset$, we get:
If $GC=D$ for all $G \lessdot D$ then $C=D$.
As another special case, when $F=D \in \C$, the second condition holds
vacuously and we get $DC=D$.

\item[(P3)]

If $FC=D$ and $C_1-C_2-\ldots-C_n=D$ is any weak $C$ gallery from $C_1$ to
$D$ then $FC_1=D$.

\end{itemize}
Axiom $(P1)(ii)$ says that $F=\emptyset$ acts as the identity on $\C$.
Also note that it is a special case of axiom $(P1)(iii)$ obtained by
setting $F=\emptyset$ and $C=D$.
In Figure~\ref{gallery}, if $F$ is such that $FC=D$ then axiom $(P3)$
says that $FC_1=D$.

\subsection*{Restriction axioms $(R)$}
For every $C \in \C$, there is a restriction map $R_C : \C \rightarrow
\F$ that satisfies 
\begin{itemize}
\item[(R1)] 

For any $C,D \in \C$, we have $R_C(C) = \emptyset$ and $R_C(D) \leq D$.

\item[(R2)] 

For any $C \in \C$, we have 
$\F = \sqcup_{D \in \C} \{ F \in \F \mid R_C(D) \leq F \leq D \}$.

\item[(R3)] 
If $R_C(C_1) \leq C_2, R_C(C_2) \leq C_3, \ldots, R_C(C_{n-1}) \leq C_n = D$
then $R_{C_1}(D) \leq R_C(D) \leq D$. 

\end{itemize}
The symbol $\sqcup$ in axiom $(R2)$ stands for disjoint union.
Also for simplicity of
notation and suggestiveness, we will write the set of inequalities in
the ``if part'' of axiom $(R3)$ as
$C_1\ \frac{R_C(C_1)}{} \ C_2\ \frac{R_C(C_2)}{}\ldots\frac{R_C(C_{n-1})}{} \ C_n$.
We will refer to this diagram as a sequence of $C$ inequalities.
In proving the equivalence of the axioms, we will see that a sequence
of $C$ inequalities 
is indeed a weak $C$ gallery (and extremal in a certain sense). 
The ambiguity in notation will then disappear. 

\subsection*{Shelling axioms $(S)$}
For every $C \in \C$, there is a partial order $\leq_C$ on
$\C$ that satisfies 

\begin{itemize}
\item[(S1)]

For any $F \in \F$, the partial order $\leq_C$ restricted to
$\C_{\geq F}$ has a unique minimal element.
For the empty face $F=\emptyset$, this unique minimal element is $C$ itself.

\item[(S2)]

Every linear extension of $\leq_C$ is a shelling of $\Delta$.

\item[(S$2^{\prime}$)]
There exists a linear extension of $\leq_C$ that is a shelling of $\Delta$.

\item[(S3)]

The partial orders are compatible in the sense that 
if $D \leq_C D_1 \leq_C D_2$ then $D_1 \leq_D D_2$.
\end{itemize}

\medskip
\noindent
Axiom $(S2)$ clearly implies axiom
$(S2^{\prime})$.
In the course of proving Theorem~\ref{t:axiom},
we will see that we can drop
axiom $(S2)$ altogether and replace it by the weaker axiom
$(S2^{\prime})$.
Hence the two axioms can be used interchangeably.

\subsection{Connection between the axioms} \label{subs:ca}

\begin{figure}[hbt]
\centering
\begin{tabular}{c@{\qquad}c}
\mbox{\epsfig{file=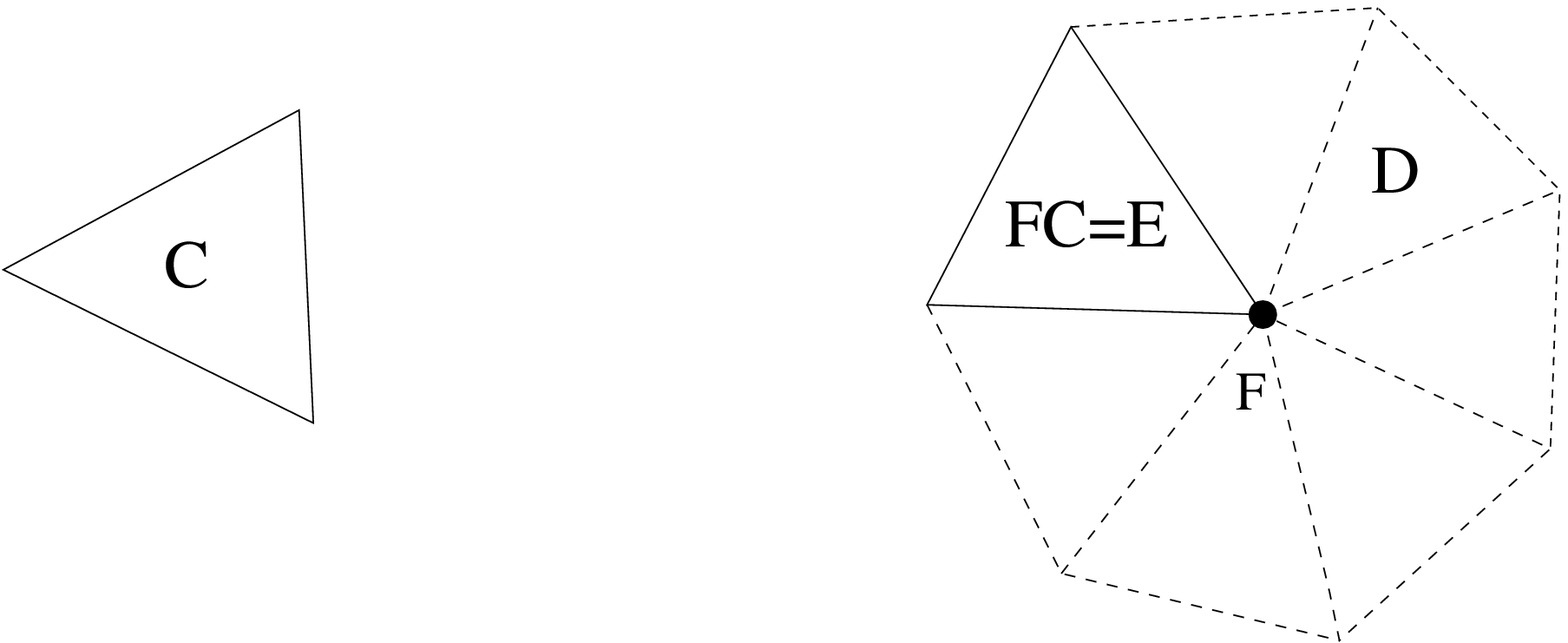,width = 8cm}}
\end{tabular}
\caption{$E$ occurs before $D$ in the shelling associated
to $C$.} 
\label{connect}
\end{figure}

We first explain the basic idea.
Projection maps on $\D$ give rise to many shellings of $\D$.
In fact we get a partial order (shelling) $\leq_C$ for every $C \in
\C$ as follows.
We say $E \leq_C D$ if there is a common face $F$ of the chambers
$D$ and $E$ such that
$FC=E$. This is illustrated in Figure~\ref{connect}.
To be technically correct, the partial order $\leq_C$ is the transitive
closure of the above relation.

\noindent
To put in words, among all chambers that contain $F$, the chamber that
is smallest in the partial order $\leq_C$ is the chamber $FC$.
And this is true for every $F \in \F$.
The restriction map $R_C$ associated to $\leq_C$ can be defined as
follows.
The face $R_C(D)$ is the smallest face $F$ of
$D$ such that $FC=D$.

We now record the above idea a little more formally. It is this formal
connection between the axioms that we will use to prove
Theorem~\ref{t:axiom}. Going from

\vspace{.2cm}
\noindent
$(P)$ to $(R)$:
For every $C \in \C$, we define 
a map $R_C : \C \rightarrow \F$.
For $D \in \C$, we let $R_C(D)$ be the face spanned by those vertices
$v$ of $D$ which satisfy $(D \setminus v) C \not= D$.

\vspace{.2cm}
\noindent
$(R)$ to $(S)$:
For every $C \in \C$, we define a partial order $\leq_C$ to be the
transitive closure of the relation:
$E \leq_C D$ if $R_C(E) \leq D$.
(Compare with equation~\eqref{e:con}.)

\vspace{.2cm}
\noindent
$(S)$ to $(P)$:
For $F \in \F$ and $C \in \C$, define $FC$ to be the unique minimal element
in $\C_{\geq F}$ with respect to the partial order $\leq_C$.
Here we used axiom $(S1)$.

\subsection{A shellable complex with a transitive group action} \label{e:nmetric}
In Section~\ref{s:proj} we discussed some examples from the viewpoint
of projection maps. Now we consider an example that is more natural
from the point of view of shellability.
It remains to be seen whether this approach 
can be formalised to get more examples.

\begin{figure}[hbt]
\centering
\begin{tabular}{c@{\qquad}c}
\mbox{\epsfig{file=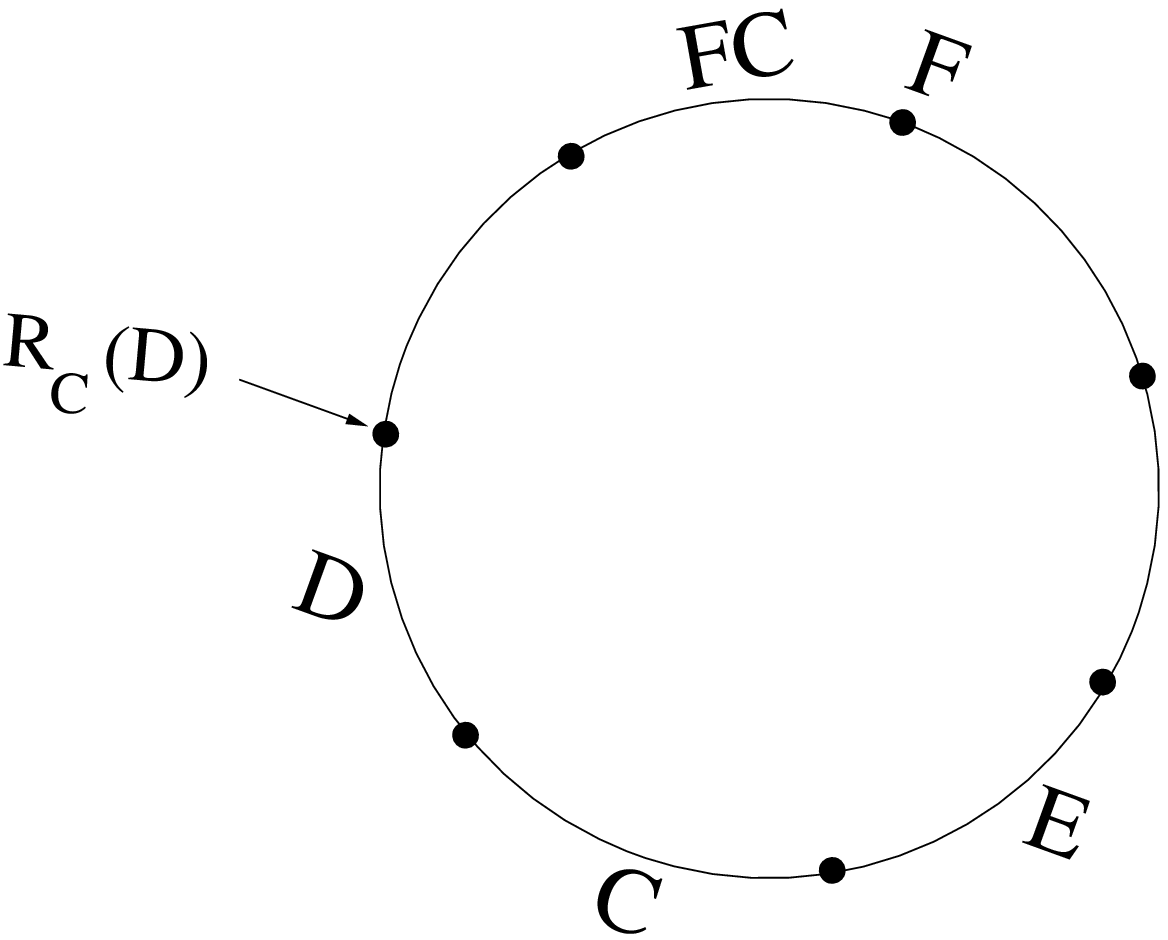,width = 4cm}}
\end{tabular}
\caption{A non-metrical example.}
\label{circle}
\end{figure}

Let $\D$ be a simplicial complex of rank $2$ that triangulates a
circle; see Figure~\ref{circle}. 
For every $C \in \C$, we define the partial order $\leq_C$ to be
the total order on $\C$ that shells $\D$ in the clockwise direction
starting at $C$. The shelling axioms are easily checked.
Note that there is a transitive action of $\Z / n\Z$ on $\D$,
where $n$ is the number of edges in $\D$. 
And this action is in some sense compatible with the partial orders $\leq_C$.

Using the connections between the three sets of axioms sketched above, we now
describe the restriction and projection maps. 
Note that $R_C(C)=\emptyset$ and $R_C(E)=E$, where 
$E$ is the chamber adjacent to $C$ in the anticlockwise direction.
For any other chamber $D$, $R_C(D)$ is the vertex of $D$ that is
further from $C$ in the clockwise direction.
For the projection maps, the only non-trivial
case is when 
$F$ is a vertex. If $F \leq C$ then $FC=C$.
If not then we define $FC$ to be the projection of $C$ on $F$ in the
clockwise direction.
Note that $FC$ is not necessarily the closest chamber to $C$ that
contains $F$. In this sense $\D$ is a non-metrical example. 
It is a good exercise to directly check the projection and restriction axioms.

\begin{remark}
In the examples that we gave in Section~\ref{s:proj}, 
we always had a map $\F \times \F \rightarrow \F$.
In our axiomatic setting, we may define such a map by using the
projection maps $\C \rightarrow \C$.
For example, we can set
$FG = \cap_{G<D,D \in \C} FD$ for $F,G \in \F$.
However it is not clear what this would imply. For instance, we may ask
whether we always get a LRB (non-associative included).
In the example above, 
if we try to extend the projection maps to a product $\F \times \F
\rightarrow \F$, then $\D$ is at best a non-associative LRB. 
\end{remark}

\subsection{Main result}
We now prove the main result of this paper.

\begin{theorem} \label{t:axiom}
Let $\Delta$ be a finite pure simplicial complex. Then
$\Delta$ satisfies $(P)$ $\Leftrightarrow$ $\Delta$ satisfies $(R)$ $\Leftrightarrow$ $\Delta$ satisfies $(S)$.
\end{theorem}
\begin{proof} 
We show $(P) \Rightarrow (R) \Rightarrow (S) \Rightarrow (P)$ 
using the connections between the axioms that we have already
outlined. The proof is fairly routine.
There is no particular reason why we choose this circle of implications.
For instance, as an exercise,
you may try to show $(P) \Rightarrow (S)$ directly.

\subsection*{$(P) \Rightarrow (R)$}
We verify axioms $(R1)$, $(R2)$ and $(R3)$.

\medskip
\noindent
$(R1)$. This is immediate from the definition of $R_C$ and axiom $(P1)(ii)$.

\medskip
\noindent
$(R2)$. Let $G$ be any facet of $D$. Then by the definition of $R_C$,
we have $G C = D \Leftrightarrow R_C(D) \leq G$.  Now by axiom $(P2)$,
we get $R_C(D) C = D$. In fact $R_C(D)$ is the unique smallest face
of $D$ with this property. This follows from axiom $(P1)(iii)$ and the
above ``if and only of'' statement. Hence $\{ F \in \F \mid R_C(D)
\leq F \leq D \} = \{ F \in \F \mid F C = D \}$.  Axiom $(R2)$ is
immediate from this description.

\medskip
\noindent
$(R3)$. 
A sequence of $C$
inequalities  
$C_1\ \frac{R_C(C_1)}{} \ C_2\ \frac{R_C(C_2)}{}\ldots\frac{R_C(C_{n-1})}{} \ C_n=D$
is in fact a weak $C$ gallery.
This is due to the fact that $R_C(C_i) C = C_i$. 
Also $R_C(D)$ is such that $R_C(D) C = D$.
Hence applying axiom $(P3)$ to the above weak $C$ gallery with $F=R_C(D)$,
we get $R_C(D) C_1 = D$.  
Since $R_{C_1}(D)$ is the smallest face $F$ such that $F C_1 = D$, we
get 
$R_{C_1}(D) \leq R_C(D) \leq D$.

\subsection*{$(R) \Rightarrow (S)$}
We first show that $\leq_C$ is a partial order.
For this we will use all three restriction axioms.
The only thing to check is antisymmetry.
Expanding out the definition of $\leq_C$ that we have given,
$C_1 \leq_C C_n$ if we have a sequence of $C$ inequalities
$C_1\ \frac{R_C(C_1)}{} \ C_2\ \frac{R_C(C_2)}{}\ldots\frac{R_C(C_{n-1})}{} \ C_n$.
To check
antisymmetry, we show that a (non-trivial) sequence as above cannot 
close in on itself.
Suppose not, then we have a minimal circular sequence, 
$${\textstyle C_1\ \frac{R_C(C_1)}{} \ C_2\ \frac{R_C(C_2)}{}\ldots\frac{R_C(C_{n-1})}{} \ C_n\ \frac{R_C(C_{n})}{} \ C_{n+1}=C_1,}$$
which is shown schematically in Figure~\ref{contradict}.
Then by axiom $(R3)$, we have $R_{C_1}(C_n) \leq R_C(C_n) \leq C_n$.
However by axiom $(R1)$ we also have $\emptyset = R_{C_1}(C_1) \leq R_C(C_n) \leq
C_1$.
Since $C_1 \not= C_n$, this contradicts axiom $(R2)$.
This shows that $\leq_C$ is a partial order.

\begin{figure}[hbt]
\centering
\begin{tabular}{c@{\qquad}c}
\mbox{\epsfig{file=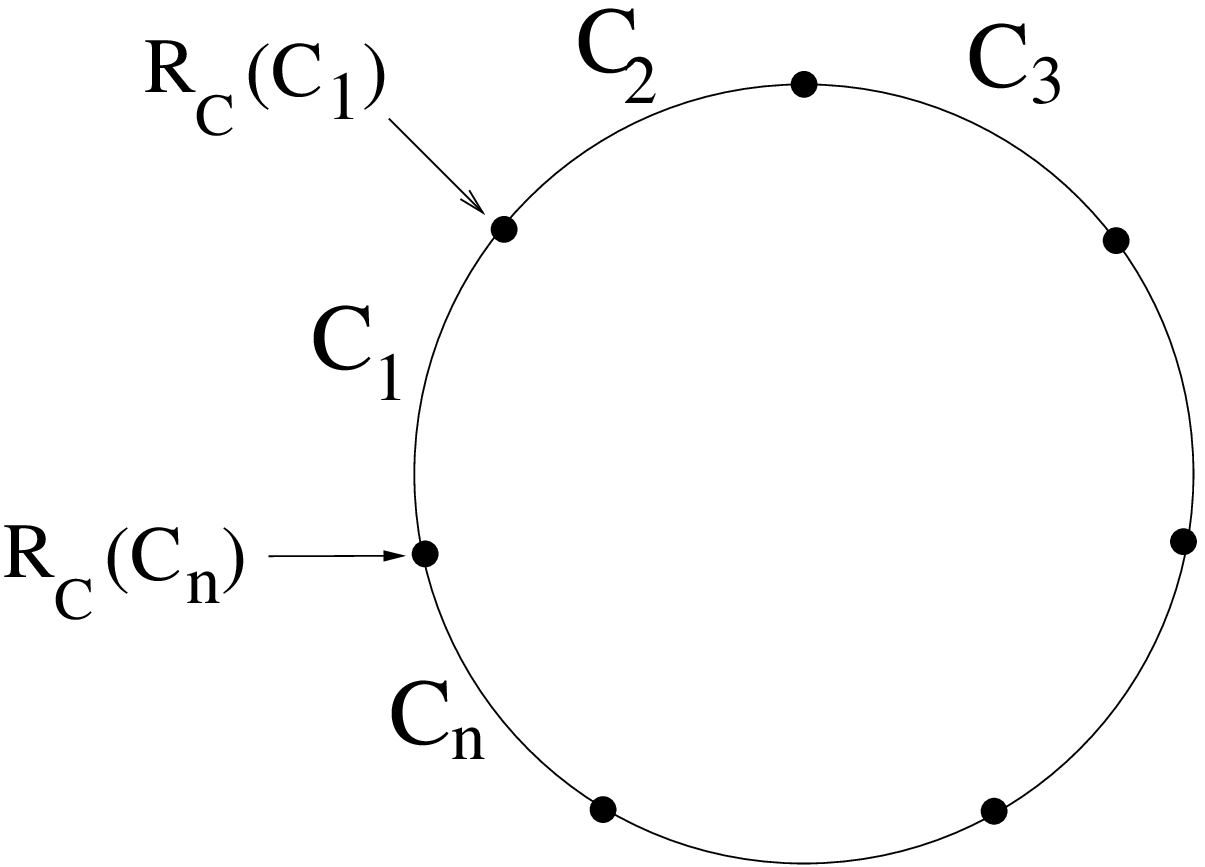,width = 4cm}}
\end{tabular}
\caption{A situation that can never occur.}
\label{contradict}
\end{figure}

\begin{remark}
In terms of projection maps, $R_C(C_n) C_1$ is forced to be both $C_1$
and $C_n$, which shows why the situation in Figure~\ref{contradict}
never occurs.
\end{remark}

\noindent
Now we verify axioms $(S1)$, $(S2)$ and $(S3)$.

\medskip
\noindent
$(S1)$. By axiom $(R2)$, for any $F \in \F$, there is a unique chamber
$D$ such that  $R_C(D) \leq F \leq D$.
For any other chamber $E$ containing $F$, we have $R_C(D) \leq
E$. Hence by definition of $\leq_C$, we have $D \leq_C E$. 
This shows that $D$ is the
unique minimal element for $\leq_C$ when restricted to 
$\C_{\geq F}$.

\medskip
\noindent
$(S2)$. Let $\leq_S$ be any linear extension of $\leq_C$. 
Since $R_C(C) = \emptyset$, we know that $C$ is the unique smallest element
in the partial order $\leq_C$. To show that $\leq_S$
is a shelling, it is enough to prove the following.\\
Claim. For any $D \in \C$ such that $D \not= C$, 
we have $D \cap (\cup_{E<_S D} E) = \cup_{v
\in R_C(D)} D \setminus v$.
Note that the RHS is a non-empty union of certain facets of $D$ as
required in the shelling condition.\\
Proof of the claim.$(\supseteq)$
Since $R_C(D) \nsubseteq D \setminus v$ for $v \in R_C(D)$, 
we see by axiom $(R2)$ that $R_C(E) \leq D \setminus v \leq E$ for some $E \not=D$. 
Now since $R_C(E) \leq D$, we have $E <_S D$.
Thus $E$ is a chamber such that $D \setminus v \leq E$ and $E <_S D$ as
required.\\
$(\subseteq)$
We show that if $F$ is a face of $D$ such that $F \nsubseteq RHS$ then
$F \nsubseteq LHS$.
If $F \nsubseteq RHS$ then $R_C(D) \leq F$. Hence if $F \leq E$ for
any $E \in \C$ then $R_C(D) \leq E$, that is, $D \leq_S E$. Hence
$F \nsubseteq LHS$.

\medskip
\noindent
$(S3)$. If $D \leq_C D_1 \leq_C D_2$ then we have a sequence of $C$
inequalities 
$D\ \frac{R_C(D)}{}\ldots-D_1\ \frac{R_C(D_1)}{}\ldots-D_2$.
Then axiom $(R3)$ implies that $D_1\ \frac{R_D(D_1)}{}\ldots-D_2$
is a sequence of $D$ inequalities. Note that we have replaced $C$ by
$D$. So $D_1 \leq_D D_2$.

\subsection*{$(S) \Rightarrow (P)$}
Now we verify axioms $(P1)$, $(P2)$ and $(P3)$.
For this we will use the weaker axiom $(S2^{\prime})$ instead of
$(S2)$.

\medskip
\noindent
$(P1)$. (i) This follows directly from the definition of the product.\\
(iii) Here we assume that $FC=D$ and $F \leq G \leq D$. The first
assumption says that $D$ is the unique minimal element in 
the partial order $\leq_C$ restricted to
$\C_{\geq F}$. 
Since $F \leq G \leq D$, it follows that $D$ is also the unique
minimal element in the partial order $\leq_C$ restricted to $\C_{\geq G}$. 
Hence $GC=D$.
As mentioned before, axiom $(P1)(ii)$ is a special case of $(P1)(iii)$.

\medskip
\noindent
$(P2)$. Using axiom $(S2')$, we choose a linear extension $\leq_S$ of $\leq_C$
which is a 
shelling of $\Delta$. 
By Lemma~\ref{l:sh}, this linear order when
restricted to $\C_{\geq F}$ gives a shelling of $\Delta_{\geq F}$, and
similarly of $\D_{\geq G}$ for all $F \leq G \lessdot D$.
Our assumption $GC=D$ for all $F \leq G \lessdot D$ 
says that for all such $G$, 
the shelling order restricted to $\C_{\geq G}$ has $D$ as the
first element. 
If we assume that $FC \not = D$ then $D$ is not the first
element in the shelling order restricted to $\C_{\geq F}$.
Hence $D \cap (\cup_{E<_S D, E \in
\C_{\geq F}} E)$ is a non-empty union of certain facets of $D$. 
This implies that there is a facet $G \lessdot D$ for which the
shelling of $\D_{\geq G}$ did not start with $D$. 
This is a contradiction. Therefore $FC=D$.

\medskip
\noindent
$(P3)$. If $D$ and $E$ are weakly $C$ adjacent then $D \leq_C
E$. Hence if $C_1-C_2-\ldots-C_n=D$ is a weak $C$ gallery then $C_1 \leq_C D$.
Let $F$ be such that $FC = D$. Now if $E \in \C_{\geq F}$ then $D
\leq_C E$. Therefore $C_1 \leq_C D \leq_C E$. 
Applying axiom $(S3)$, we get $D \leq_{C_1} E$ for all $E \in \C_{\geq
F}$. Hence $FC_1=D$.

\end{proof}

From now on, if $\D$ is a pure simplicial complex that satisfies $(P)$,
or equivalently, $(R)$ or $(S)$ then we will simply say that $\D$
satisfies our axioms.

\section{An almost example} \label{s:LRB}

The purpose of this section is to identify a potential source of
examples that satisfy our axioms. 
In Section~\ref{s:axiom} the poset $\F$ was the face lattice of a simplicial
complex $\D$. 
It is natural to consider a more general situation,
where instead of $\D$, we just have a poset $\F$ with the
maximal elements playing the role of $\C$.
The projection axioms still make sense though the role they now play
is not exactly clear. 
An interesting case to consider is when $\F$ is a LRB.
We will show that a LRB satisfies axioms $(P1)$ and $(P3)$.
However it can easily violate axiom $(P2)$.
We first discuss this situation.

\begin{figure}[hbt]
\centering
\begin{tabular}{c@{\qquad}c}
\mbox{\epsfig{file=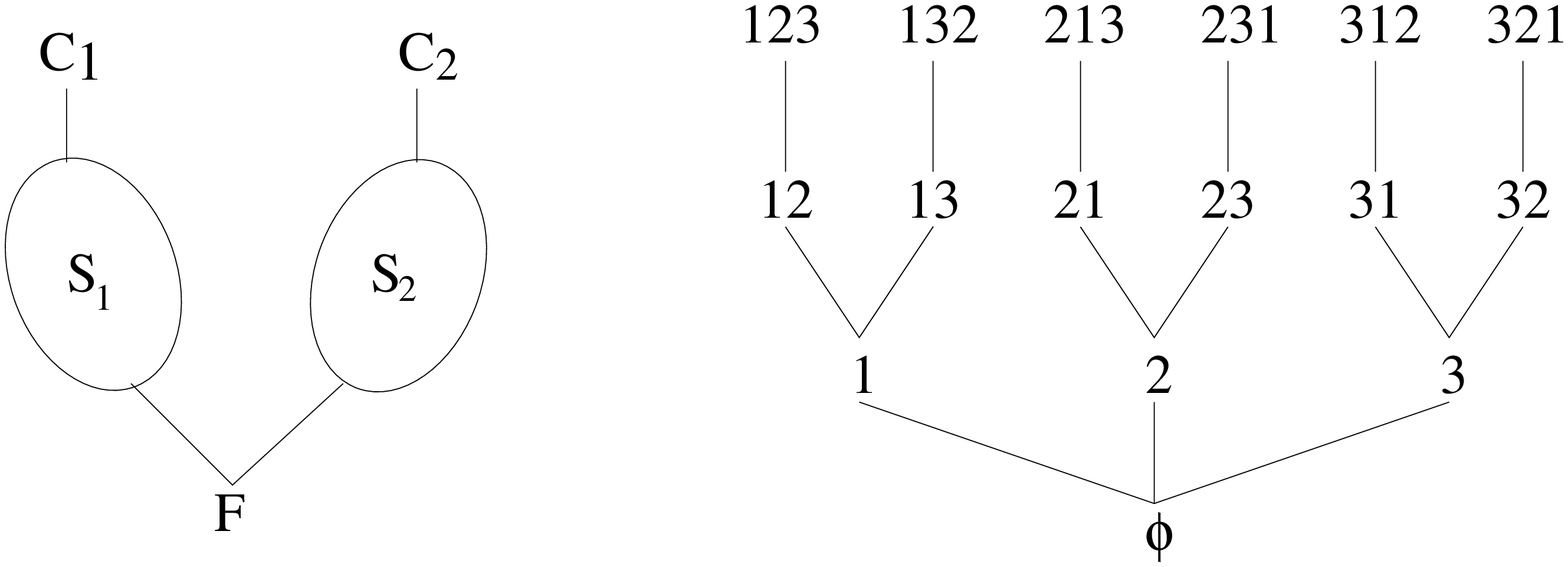,height=3cm, width = 10cm}}
\end{tabular}
\caption{}
\label{lrb}
\end{figure}

\subsection*{On why a LRB violates axiom $(P2)$}
To see what goes wrong, consider the following example.
Let $S_1$ and $S_2$ be any two LRBs.
Consider the LRB $S$ shown in Figure~\ref{lrb} (on the left) with $F_1 F_2 = C_1$
and $F_2 F_1 = C_2$  for $F_1 \in S_1$ and $F_2 \in S_2$.
We first check that $S$ is indeed a LRB.
Then by definition $G C_1 = C_2$ for $F < G \lessdot C_2$.
However $F C_1 = C_1 \not= C_2$.
Hence $S$ violates axiom $(P2)$.

To give a concrete example that illustrates the same problem, consider
the \emph{free} LRB with identity 
on $n$ generators, denoted $F_n$.
The elements of $F_n$ are sequences
$x=(x_1,\dots,x_l)$ of distinct elements of the set
$[n]=\{1,\dots,n\}$, $0\le l\le n$.  We multiply two such sequences by
\[
(x_1,\dots,x_l)(y_1,\dots,y_m)=(x_1,\dots,x_l,y_1,\dots,y_m)\sphat\,,
\]
where the hat means ``delete any element that has occurred earlier''.
For example,
\[
(2\,1)(3\,5\,4\,1\,6)=(2\,1\,3\,5\,4\,6).
\]
Figure~\ref{lrb} (on the right) shows the Hasse diagram of $F_3$, the free
LRB on $3$ generators. 
Check that $F_n$ violates axiom $(P2)$.
Also note that every chamber has only one facet. 

These examples show that a LRB
is not a geometric object in any sense of the term.
A natural question that comes up is whether a simplicial LRB always
satisfies axiom $(P2)$. We do not know the answer.
For the LRBs that we considered in Section~\ref{s:proj},
the answer is yes.
In other words, the LRB associated to any hyperplane arrangement
satisfies axiom $(P2)$ (and hence the same holds for flag (order) complex
$\D(M)$ of a distributive lattice $M$). 
It is an easy and illuminating exercise to check this directly using the
definition of the product involving sign sequences.

\subsection*{A LRB satisfies axioms $(P1)$ and $(P3)$}

Axiom $(P1)$ follows directly from the definition of a LRB and axiom
$(P3)$ holds by the following proposition. 
It is a good exercise to check axiom $(P3)$ 
directly for the free LRB.

\begin{proposition}
Let $\F$ be a poset with an associative product $\F \times \F
\rightarrow \F$. Let $\C$, the set of maximal elements of $\F$, 
be a left ideal in $\F$.
Also assume that $\F$ satisfies axiom $(P1)$ with the additional
property that $F \leq FF^{\prime}$ for any $F,F^{\prime} \in \F$.
Then $\F$ satisfies axiom $(P3)$.
\end{proposition}
\begin{proof}
Though our setting is somewhat abstract, it would help to keep in mind
the picture in Figure~\ref{gallery}.
Let 
${\textstyle C_1\ \frac{F_1}{} \ C_2\ \frac{F_2}{}\ldots\frac{F_{n-1}}{} \ C_n=D}$
be a weak $C$ gallery; that is, $F_i C = C_i$ and $F_i \leq
C_i,C_{i+1}$.
We want to show that if $F C = D$ then $F C_1 = D$.
For notational consistency, we put $F=F_n$.
Let $\G = \{ H \in \F \mid F_n \leq H \leq C_n = F_n C \}$.
We need to show that $F_n C_1 = C_n$.

We first note three easy consequences of parts $(i)$, $(ii)$ and $(iii)$
respectively of axiom $(P1)$.

$(1)$. If $F_n \leq H$ and $H E = C_n$ for some $E \in \C$ then $H \in
   \G$.

$(2)$. For $2 \leq i \leq n$, we have $F_{i-1} F_i C = F_i C$. This is
   because $F_{i-1} \leq  C_i = F_i C$.

$(3)$. If $H \in \G$ then $H C = C_n$.

\medskip
\noindent
We now claim that $F_n F_i \in \G$ for $1 \leq i \leq n$.\\
Proof of the claim. We do a reverse induction on the index $i$.
Clearly $F_n F_n \in \G$.
To do one more step, note that $F_{n} F_{n-1} \in \G$. 
This is because $F_n \leq F_n F_{n-1}$ by our additional assumption
and $F_n F_{n-1} \leq C_n$ follows from $F_{n-1},F_n \leq C_n$ and
axiom $(P1)$.
Now assume that $F_n F_i \in \G$ for some $i$. We show that $F_n
F_{i-1} \in \G$. 

Using $(3)$ we get $(F_n F_i) C = C_n$.
Rewriting this as $F_n (F_i C) = C_n$ and using $(2)$, we get
$F_n (F_{i-1} F_i C) = C_n$ which we write as
$(F_n F_{i-1}) F_i C = C_n$.
Note that $F_n \leq F_n F_{i-1}$.
Hence we apply $(1)$ with $H = F_n F_{i-1}$ and $E = F_i C$ to conclude
that $F_n F_{i-1} \in \G$. This completes the induction step and the
claim is proved.

From the claim we get $F_n F_{1} \in \G$. 
Using $(3)$ we have $(F_n F_1) C = C_n$ which is same as
$F_n (F_1 C) = F_n C_1 = C_n$.
This is exactly what we wanted to show.
\end{proof}

\noindent
The discussion in this section shows that the LRBs in Section~\ref{s:proj}
satisfy the projection axioms and hence are examples of our theory.
In the next section, we again consider these examples 
(including the non-associative ones) but from a
more intuitive and geometric perspective.

\section{The metric setup} \label{s:metric}

In this section we show that if a chamber complex $\D$
satisfies the gate property then it
satisfies our axioms.
This is essentially an axiomatic restatement of a result of Scharlau;
see the proposition in~\cite[Section 3]{\scharlau}.
(He states it in terms of shellings.
Also he does not consider the compatibility axiom.)
Furthermore in this special situation, the maps of
Section~\ref{s:axiom} have a geometric meaning. 
To make this clear we first need some definitions.

\subsection{Some definitions}
Let $\D$ be a pure simplicial complex. 
The maximal simplices are called \emph{chambers}.
We say two chambers are \emph{adjacent} if they have a common codimension $1$
face. A \emph{gallery} is a sequence of chambers such that consecutive
chambers are adjacent. We say that $\D$ is \emph{gallery connected} if for
any two chambers $C$ and $D$, there is a gallery from $C$ to $D$.
For any $C,D \in \C$, we then define the gallery distance $d(C,D)$ to
be the minimal length of a 
gallery connecting $C$ and $D$. And any gallery which achieves this
minimum is called a \emph{geodesic gallery} from $C$ to $D$.
Another relevant metric notion is the following.

\medskip
\noindent
\emph{Gate Property.} For any face $F \in \F$ and chamber $C \in \C$,
there exists a chamber $D \in \C_{\geq F}$ such that
$d(C,D) \leq d(C,E)$ for any $E \in \C_{\geq F}$.
Furthermore
$d(C,E) = d(C,D) + d(D,E)$.

\medskip
\noindent
The gate property implies that $\D$ is \emph{strongly connected}; that is,
$\C_{\geq F}$ is gallery connected for all $F \in \F$.
In fact it implies that $\C_{\geq F}$ is a \emph{convex} subset of $\C$; 
that is, if $D$ and $E$ are any two chambers in $\C_{\geq F}$ then any
geodesic gallery from $D$ to $E$ lies entirely in $\C_{\geq F}$.

\begin{figure}[hbt]
\centering
\begin{tabular}{c@{\qquad}c}
\mbox{\epsfig{file=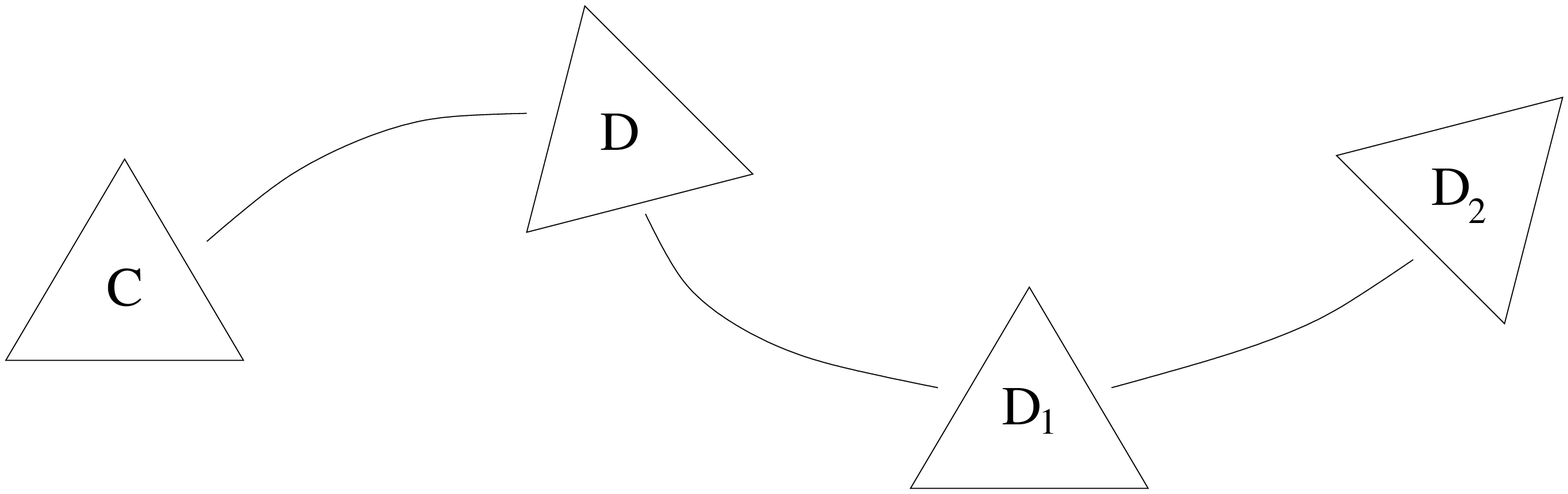,height=2.3cm,width = 8cm}}
\end{tabular}
\caption{A geodesic gallery that illustrates shelling compatibility.}
\label{compatible}
\end{figure}

\subsection{Geometric descriptions of the relevant maps}

We now use the above metric notions to define the maps that played a
role in the axiomatic setup of Section~\ref{s:axiom}.
Let $\D$ be a \emph{chamber complex}. 
It is by definition a gallery connected
pure simplicial complex.
For every $C \in \C$, we define a partial order $\leq_C$.
We say that $D \leq_C E$ if there is a geodesic gallery from $C$ to
$E$ that passes through $D$.
With this definition of the partial orders $\leq_C$, 
we analyse the shelling axioms.
We first note that axiom $(S3)$ is a consequence of our definition.
If $D \leq_C D_1 \leq_C D_2$ then we have a geodesic gallery 
$C-\ldots-D-\ldots-D_1-\ldots-D_2$ as shown in Figure~\ref{compatible}.
This restricts to a geodesic gallery
$D-\ldots-D_1-\ldots-D_2$, which implies $D_1 \leq_D D_2$.
Now we look at axiom $(S1)$. It says that for any $C \in \C$ and $E
\in \C_{\geq F}$ there is a unique chamber $D \geq F$ such that $D \leq_C
E$; that is, such that there is a geodesic gallery $C-\ldots-D-\ldots-E$.
This is equivalent to saying that $\D$ has the gate property.
Next we claim that a chamber complex $\D$ with
the partial orders $\leq_C$ as defined above 
satisfies the shelling axioms $(S)$ if and
only if it has the gate property.
From what we have so far, we only need to show that the gate property implies
axiom $(S2)$. This can be checked directly. However we will prove this using
the projection axioms.

\begin{remark}
If $\D = \Sigma$ is the Coxeter complex associated to a Coxeter group $W$
then $\leq_C$ coincides with the weak Bruhat order on $W$, after we
choose $C$ as the fundamental chamber and identify chambers of
$\Sigma$ with elements of $W$.
\end{remark}

Now we assume that $\D$ has the gate
property. We show that this implies the projection axioms. 
Using the gate property, we define $F C$ to be the chamber containing
$F$ that is closest to $C$. Axiom $(P1)$ follows. The gate property further
says that a weak $C$ gallery $C_1-C_2-\ldots-C_n$ 
can be extended to a geodesic gallery $C-\ldots-C_1-\ldots-C_2-\ldots-C_n$. 
Axiom $(P3)$ now follows. To show 
axiom $(P2)$, we proceed by contradiction. 
Let $FC = E \not= D$.
Choose a geodesic gallery
$C-\ldots-E-\ldots-D$.
Let $G$ be the facet of $D$ that is crossed by this gallery in the
final step.
Then $GC \not= D$, which is a contradiction.
To summarise, $\D$ satisfies the projection axioms $(P)$ 
(with the closest chamber
definition) if and only if it satisfies the gate property.

To complete the story we now describe the restriction maps. 
Recall from Section~\ref{subs:ca} that $R_C(D)$ is the face of
the chamber $D$ spanned by vertices $v$ such that $(D \setminus
v)C\not=D$. It can be equivalently described as the face of the chamber
$D$ spanned by vertices $v$ which have 
the following property.

\medskip
\noindent
There is a minimal gallery from $C$ to $D$ that passes through the
facet $D \setminus v$ in the final step.

\medskip
\noindent
We refer to the second description as $\Des(C,D)$. 
The equivalence of the two descriptions says that $R_C(D)=\Des(C,D)$.
The terminology $\Des(C,D)$ is again motivated 
by the theory of Coxeter groups. 
$\Des(C,D)$ is the face of $D$ spanned by ``the descent set of $D$ with
respect to $C$''. For an explanation, see the appendix by Tits to Solomon's
paper~\cite{\solomon} or the more elaborate exposition in 
\cite[Section 9]{\ken}.

We have proved the following theorem.

\begin{theorem} \label{t:maxiom}
Let $\D$ be a chamber complex.
Also let the partial orders $\leq_C$ , the restriction maps $R_C$
and the projection maps be as defined in this section. Then $\D$
satisfies our axioms
if and only if it
satisfies the gate property. 
\end{theorem}

\subsection{Examples of Section~\ref{s:proj} revisited}
Hyperplane arrangements, buildings and modular lattices are all 
examples of chamber complexes that
satisfy the gate property.
Hence by Theorem~\ref{t:maxiom} 
we have established that they indeed satisfy our axioms.

For oriented matroids (in particular, central hyperplane arrangements), 
the gate property is due to Bj{\"o}rner and Ziegler and
for buildings it is due to Tits~\cite{\tits}. 
The gate property for modular lattices 
(in particular, distributive lattices) was proved by Abels, see
~\cite[Proposition~2.9]{\abels}.
Alternatively, for distributive lattices,
we can deduce the gate property by combining the following facts.
If $\D$ satisfies the gate property then so does any subcomplex
$\D_c$ spanned by a convex subset of chambers.
The flag (order) complex $\D(M)$ of a distributive lattice $M$
corresponds to a convex set of chambers in $\D(\mathcal{B}_{n})$, 
the complex of the braid arrangement.

\medskip
\noindent
\emph{More on Example 1}.
The metric notions that we described in this section can be made
very explicit for this example.
We prefer to restrict to the simplicial case
though all statements except the last make sense and hold in general.
A minimal gallery from $C$ to $D$ crosses exactly those hyperplanes
that separate $C$ and $D$.
Consequently the
gallery distance $d(C,D)$ is given by the number of hyperplanes
separating $C$ and $D$.
We have $D \leq_C E$ if and only if the hyperplanes that separate $C$
and $D$ also separate $C$ and $E$.
The chamber $FC$  is the unique chamber containing $F$ for which no hyperplane
passing through $F$ separates $C$ from $FC$.
Also 
$R_C(D)$ is the face of $D$ whose support is the intersection of those
walls of $D$ that do not separate $D$ from $C$.
Rigorous proofs of these unjustified facts 
can be found in any of the references cited earlier or
you may accept them as ``intuitively obvious''.

\section{Type selected subcomplexes} \label{s:tss}

Many results of this section are well-known and in greater generality
than what we present here. For example, see~\cite[Section 8.3]{\zie}. 
Even then we give a self-contained treatment mainly because our
approach is different and uses projection maps.
Also we will use some of these ideas to motivate
the discussion in the next three sections.

A \emph{labelling} of a pure simplicial complex $\D$ by a set $I$ is a
function which 
assigns to each vertex an element of $I$, in such a way that the
vertices of every chamber are mapped bijectively onto $I$.
A labelled simplicial complex is also sometimes called ``numbered'' or
``balanced''. 

Let $\D$ be a labelled simplicial complex and let $\D_J$ be
its \emph{type selected} subcomplex consisting of faces whose label set or
\emph{type} is contained in $J$. We further assume that $\D$
satisfies our axioms.
We first study the connection between the homotopy type of $\D_J$
and the flag $h$ vector. To get this far, we do not really need the
compatibility axiom. So we fix a chamber $C \in \C$, assume the first
two axioms (in each category) for $C$ and work from there.
However in Section~\ref{ss:h_ddesc},
we make full use of the compatibility axiom to look at the dual
picture for the flag $h$ vector. The fact that this viewpoint is useful can be seen from the
example that we present later in Section~\ref{ss:example}.
We also consider the unlabelled case briefly in Section~\ref{ss:unlabelled}.

We start with the following result of Bj{\"o}rner 
that generalises Proposition~\ref{p:shell}. 

\begin{proposition}
\cite[Theorem 1.6]{\bjorner} \label{p:type}
For a labelled shellable complex $\D$,
the type selected subcomplex $\D_J$ is shellable and homotopy
equivalent to a wedge of 
$(\abs{J}-1)$-spheres. The number of these spheres is given by
$\beta_J = \abs{\{D \mid R(D) \  \text{has type} \  J\}}$, 
where $R$ is the restriction map associated with the shelling.
\end{proposition}

\begin{remark}
The proposition is proved by showing that the shelling order for $\D$
``induces'' a shelling order for $\D_J$.
This raises the following question.
Suppose we assume that $\D$
satisfies our axioms with $R_C$ as the restriction maps.
Also let $R_{J,C}$ be the maps for $\D_J$
induced by $R_C$.
Then what kind of restriction axioms do the induced maps $R_{J,C}$ 
satisfy for $\D_J$?
The axioms need to be suitably generalised because 
there is more than one restriction map for every
chamber in $\D_J$.
\end{remark}

The proposition applies to our situation because of axiom $(S2)$ which
implies that $\D$ is shellable.
The restriction map $R$ in the proposition is the map $R_C$ for us, 
where $C$ is the chamber that we have fixed.

\subsection{The restriction map and the flag $h$-vector}\label{ss:h_desc}
Recall that $R_C(D)$ is the unique smallest face $F$ of $D$
such that $FC=D$.
Next define 
$$\beta_J = \abs{\{D \mid R_C(D) \ 
\text{has type} \  J\}}.$$ 
We will show that the vector $(\beta(J))_{J\subseteq
I}$ coincides with the ``flag $h$-vector'' defined below.

First we define the \emph{flag $f$-vector} of $\D$ by setting
$f_J(\D)$ equal to the number of simplices of type~$J$.  The
\emph{flag $h$-vector} is then obtained by writing
\begin{equation} \label{e:hvect}
f_J(\D)=\sum_{K\subseteq J} h_K(\D),
\end{equation}
or, equivalently,
\begin{equation}
h_J(\D)=\sum_{K\subseteq J} (-1)^{\abs{J-K}} f_K(\D).
\end{equation}

\begin{proposition} \label{p:beta}
Let $\D$ be as above and 
let $I$ be the set of types of vertices.  Then for any $J\subseteq I$,
\[
\beta_J=h_J(\D).
\]
\end{proposition}

\begin{proof}
Let $\F_J$ be the set of simplices of type $J$.  There is a 1--1
map $\F_J\to\C$, given by $F\mapsto FC$, where $C$ is fixed.  It is
1--1 because we can recover $F$ from $FC$ 
as the face of type~$J$.  Its image is the set of chambers $D$ for
which the type of $R_C(D)$ is contained in~$J$.  Hence
\[
f_J(\D)=\sum_{K\subseteq J} \beta_K.
\]
The proposition now follows from \eqref{e:hvect}.
\end{proof}

\begin{remark}
The above argument is due to Ken Brown~\cite{\ken}. 
He states this proposition for Coxeter complexes.
Also he works with descent sets rather than restriction maps.
We know from the discussion in Section~\ref{s:metric} that in the
metric setup, $R_C(D)=\Des(C,D)$.
\end{remark}

\subsection{The dual picture}\label{ss:h_ddesc}
Now we assume that $\D$
satisfies all three axioms.
In other words, we are now allowed to vary $C \in \C$ and have the
compatibility axiom.
Note that $\beta_J$ does not depend on the choice of the chamber $C$.
For example, from Proposition~\ref{p:type}, we know that
$\beta_J$ counts the number of spheres in $\D$ and hence is
independent of the choice of $C$. Also from Proposition~\ref{p:beta},
we know that it coincides with the flag $h$ vector of $\D$.
We also see that
$h_J(\D) = \abs{\{D \mid R_C(D) \  \text{has type} \  J\}}$. 
This gives geometric interpretations of the flag $h$ vector, one for
every $C \in \C$. 

Now we make full use of the compatibility axiom by turning the picture
around. 
For the map $\F \times \C \rightarrow \C$
that maps $(F,C)$ to $FC=D$, we fix $D$ and vary $C$.
More precisely, we let
$h_J(D) = \abs{\{C \mid R_C(D) \  \text{has type} \  J\}}$. 
Similarly let
$f_J(D)$ be the number of chambers $C$ such that
$FC=D$, where $F$ is the face of $D$ of type~$J$.  
This can be rewritten as 
$f_J(D) = \abs{\{C \mid \text{The type of}\ R_C(D) \  \text{is
contained in} \  J\}}$. 
Now observe that 
\begin{equation} \label{e:dhvect}
f_J(D)=\sum_{K\subseteq J} h_K(D).
\end{equation}

\noindent
The numbers $f_J(D)$ and $h_J(D)$ do in general depend on $D$.
If we average them over all $D \in \C$ then we recover the usual flag
vectors. To see this, 
let $\sigma_J$ be the sum of all faces of type $J$. In particular,
$\sigma_I$ is the sum of all chambers.
Note that $\sigma_J \sigma_I = \sum_{D \in \C} f_J(D) D$.
By counting the total number of terms involved, we get 
$f_J(\D) = \frac{1}{\abs{\C}}\sum_{D \in \C} f_J(D)$.
This along with \eqref{e:hvect} and \eqref{e:dhvect} says that 
$h_J(\D) = \frac{1}{\abs{\C}}\sum_{D \in \C} h_J(D)$.
Note that it is not obvious from the definition that the right hand
side is an integer. 

The numbers $h_J(D)$ give us a different way 
to understand the flag $h$ vector. 
A natural question that arises is: 
When is $h_J(\D) = h_J(D)$ for all $D \in \C$?
We will address this in more detail in Section~\ref{s:comm}.

\subsection{The unlabelled case} \label{ss:unlabelled}
If $\D$ is not labelled then we work with the ordinary $f$ and $h$
vectors rather than the flag vectors.
Let $\D$ be a simplicial complex of rank $n$ that
satisfies our axioms.
As before, we will need the compatibility axiom only when we pass to the
dual picture.
We define $\beta_j = \abs{\{D \mid R_C(D) \ 
\text{has rank} \  j\}}$. 
The unlabelled version of Proposition~\ref{p:type} is as follows.

\begin{lemma}
For a shellable complex $\D$ of rank $n$,
the $(k-1)$-skeleta $\D_{k-1}$ is shellable and homotopy
equivalent to a wedge of $(k-1)$-spheres. The number of these spheres
is given by 
$\sum_{i=0}^{k} \binom{n-i-1}{k-i} \beta_i.$
\end{lemma}

The fact that shellability is inherited by $k$-skeleta appears as
Corollary 10.12 in~\cite{\bw}. Also see the references cited therein.
The formula for the number of spheres was told to me by Bj{\"o}rner.
As before, 
we can show that the vector $(\beta_j)_{0 \leq j \leq {n}}$ 
coincides with the ``$h$-vector'' defined below.

First we define the \emph{$f$-vector} of $\D$ by setting
$f_j(\D)$ equal to the number of simplices of rank~$j$.  The
$h$-vector is then obtained by writing
\begin{equation} \label{e:uhvect}
f_j(\D)=\sum_{k=0}^{j} \binom{n-k}{n-j} h_k(\D).
\end{equation}

\vanish{
or, equivalently,
\begin{equation}
h_J(\D)=\sum_{K\subseteq J} (-1)^{\abs{J-K}} f_K(\D).
\end{equation}
}

\begin{proposition} \label{p:ubeta}
Let $\D$ be as above.
Then for any $0 \leq j \leq n$,
\[
\beta_j=h_j(\D).
\]
\end{proposition}

\begin{proof}
Let $\F_j$ be the set of simplices of rank $j$.  
As before, we consider the
map $\F_j\to\C$, given by $F\mapsto FC$, where $C$ is fixed.  
Recall that for a fixed $C$ and $D$,
$\{ F \in \F \mid R_C(D)
\leq F \leq D \} = \{ F \in \F \mid F C = D \}$. 
This says that the image of the above map is the set of chambers $D$ for
which the rank of $R_C(D)$ is less than or equal to ~$j$.  
However this map is no longer 1--1.
The number of times a fixed chamber $D$ occurs in the image is
$\binom{n-k}{n-j}$. This is the number of faces of
$D$ of rank $j$ that contain $R_C(D)$ (whose rank we set equal to $k$).
This gives
\[
f_j(\D)=\sum_{k=0}^{j} \binom{n-k}{n-j} \beta_k
\]
and the proposition now follows.
\end{proof}

To get the dual picture, let
$h_j(D) = \abs{\{C \mid R_C(D) \  \text{has rank} \  j\}}$. 
Similarly let
$f_j(D) = \abs{\{(F,C) \mid FC=D,\ \text{rk} \ F=j\}}$. 
As in the proof of Proposition~\ref{p:ubeta},
$f_j(D)=\sum_{k=0}^{j} \binom{n-k}{n-j} h_k(D)$.
Now if we let $\sigma_j$ be the sum of all faces of rank $j$ then 
$\sigma_j \sigma_n = \sum_{D \in \C} f_{j}(D) D$.
This gives
$f_j(\D) = \frac{1}{\abs{\C}}\sum_{D \in \C} f_{j}(D)$
and then
$h_j(\D) = \frac{1}{\abs{\C}}\sum_{D \in \C} h_{j}(D)$.

\begin{remark}
For a labelled complex, we have 
$f_j=\sum_{\abs{J}=j} f_J, h_j=\sum_{\abs{J}=j} h_J$ and
$\sigma_j=\sum_{\abs{J}=j} \sigma_J$.
So in this case, the results for the $f$ and $h$ vector
follow from the corresponding result for the flag vectors.
For instance,
using the first two equations, one can recheck that
equation~\eqref{e:hvect} reduces to equation~\eqref{e:uhvect}.
\end{remark}

\subsection{Simplicial hyperplane arrangements} \label{ss:example}
Let $\D$ be a simplicial hyperplane arrangement.
We now give a description of the flag $h$ vector in this case. 
The dual description in Section~\ref{ss:h_ddesc} 
that involves $h_J(D)$ is easier to
visualise and more useful than the one in Section~\ref{ss:h_desc}.
Recall that in this situation, 
$R_C(D)$ is the face of $D$ whose support is the intersection of those
walls of $D$ that do not separate $D$ from $C$;
see the description at the end of Section~\ref{s:metric}.
Using this description, in the labelled case,
$h_J(D)$ is the cardinality of the set of chambers that, among the
walls of $D$, lie onside exactly those that contain
the face of type $J$ of $D$.

\begin{figure}[hbt]
\centering
\begin{tabular}{c@{\qquad}c}
\mbox{\epsfig{file=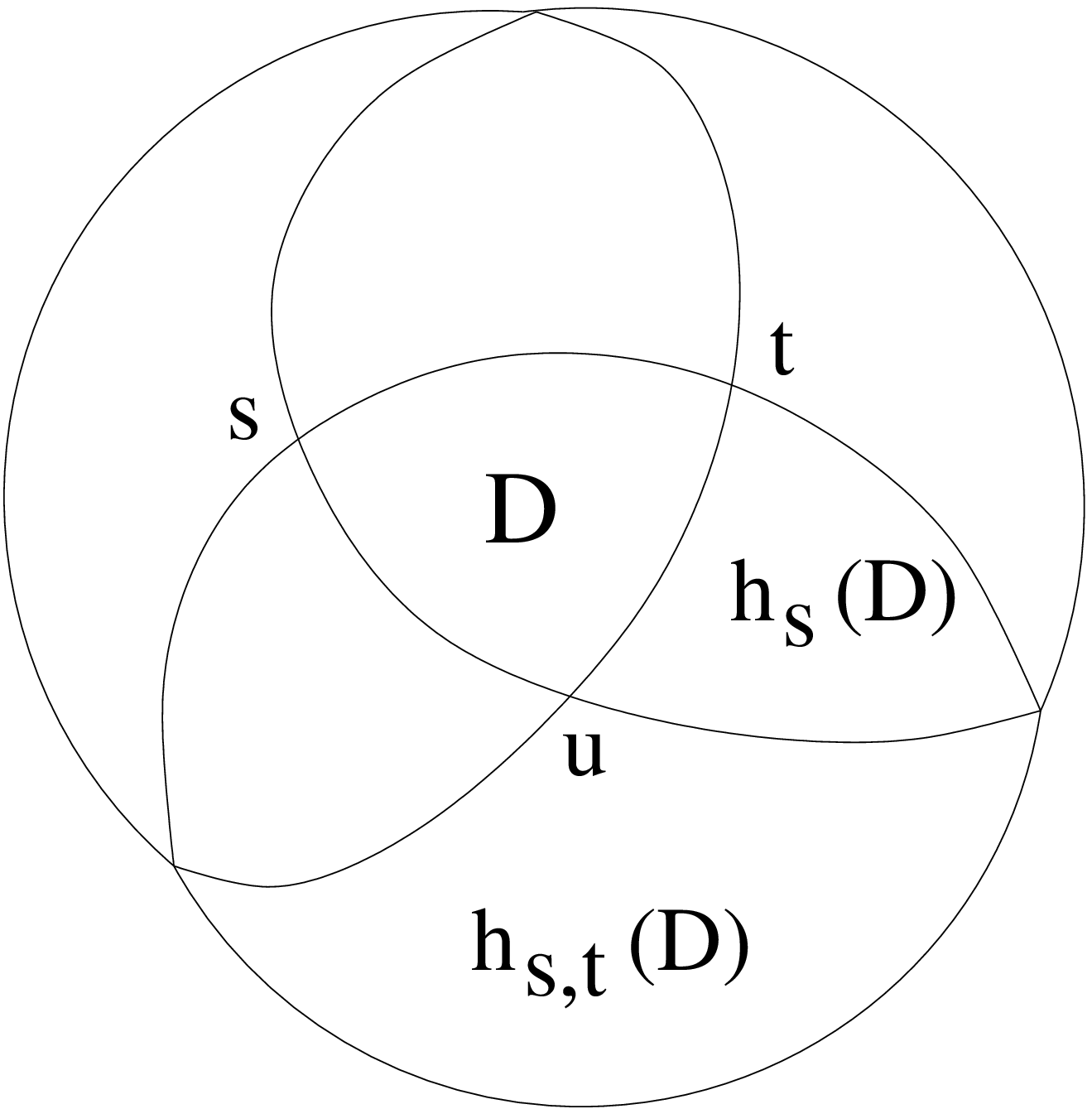,height=3cm,width = 3cm}}
\end{tabular}
\caption{The schematic picture in rank 3.}
\label{schema}
\end{figure}

\noindent
In Figure~\ref{schema} only the three walls 
of the chamber $D$ have been drawn. 
The only chamber that is not seen in the picture is the one directly
opposite to $D$. The flag $h$ vector for $D$ is represented by the 8
regions into which the three walls divide the sphere. We
have indicated the regions that correspond to $h_s(D)$ and $h_{s,t}(D)$.

In the unlabelled case,
$h_j(D)$ is the cardinality of the set of chambers that lie offside
exactly $j$ of the $n$
walls of $D$.
In Figure~\ref{schema}, for example, $h_1(D)$ counts the chambers that lie in
the three regions adjacent to $D$ along its facets.
Next we give a tiny application of this description.
Let $\D_{min}$ be the Coxeter complex of the reflection arrangement
$x_i = 0$, where $1\le i\le n$.
In fact Figure~\ref{schema} as drawn corresponds to the case $n=3$.
Observe that $h_j(D)=\binom{n}{j}$.
From our interpretation, it is clear that
for any other simplicial arrangement $\D$ of rank $n$,
we get $h_j(D) \geq \binom{n}{j}$ for all $D \in \C$ and hence
$h_j(\D) \geq \binom{n}{j}$.

Now consider the barycentric subdivision $\D^{\prime}$ of $\D$.
Then $\D^{\prime}_{min}$ corresponds to the Coxeter arrangement of type $B_n$,
namely, $x_i= \pm x_j$ and $x_i = 0$, where $1\le i<j\le n$.
However in general, $\D^{\prime}$ is not a hyperplane arrangement. 
It is just a labelled simplicial complex.
Hence we cannot directly conclude that
$h_J(\D^{\prime}) \geq h_J(\D^{\prime}_{min})$ (which is known to be true).
This raises an important question. 
If $\D$ satisfies our axioms then what can we say about its 
barycentric subdivision $\D^{\prime}$ ?
For a positive result in this direction, see \cite[Theorem 5.1]{\bjor}
which says that if $\D$ is shellable then so is $\D^{\prime}$.

\section{Commutativity issues} \label{s:comm}

In this section, we return to
the question: When is $h_J(D) = h_J(\D)$ for all $D \in \C$?
This can be rephrased as a commutativity problem.
Recall that $\sigma_J$ (resp. $\sigma_j$) is the sum of all faces of
type $J$ (resp. rank $j$). 
Also $h_J(D)$ and $h_J(\D)$ are the local and global versions of the
flag $h$ vector of $\D$.
And $h_j(D)$ and $h_j(\D)$ are the corresponding versions of the
$h$ vector of $\D$.

\begin{proposition} \label{p:equi}
Let $\D$ be a simplicial LRB (possibly non-associative)
that satisfies our axioms. Then\\
(labelled case) $h_J(\D) = h_J(D)$ for all $D \in \C$ and for all subsets $J
\subseteq I$ $\Leftrightarrow$
$\sigma_I \sigma_J = \sigma_J \sigma_I$ for all subsets $J
\subseteq I$.\\
(unlabelled case) $h_{j}(\D) = h_{j}(D)$ for all $D \in \C$ and
for all $0 \leq j \leq n$ $\Leftrightarrow$
$\sigma_n \sigma_j = \sigma_j \sigma_n$
for all $0 \leq j \leq n$.
\end{proposition}
\begin{proof}
We do only the labelled case. The unlabelled case is similar.
Suppose that
$\sigma_I \sigma_J = \sigma_J \sigma_I$ for all subsets $J
\subseteq I$.
Counting the number of times a fixed chamber $D$ occurs on both sides
of the equality, we get
$f_J(\D) = f_J(D)$, or equivalently, $h_J(\D) = h_J(D)$.
\end{proof}

\begin{remark}
In the proposition above, a concrete case to keep in mind is that of
simplicial hyperplane
arrangements. 
Also we required $\D$ to satisfy our axioms because we only defined
$h_J(D)$ in that situation.
And we required $\D$ to be a LRB so as to make sense of the product
$\sigma_I \sigma_J$.
Of course, we only needed to use that $CF=C$ for $C \in \C$ and $F \in \F$.

If $\D$ is a finite Coxeter complex or a finite (spherical) building
associated with a BN-pair then 
there is a type and product preserving group
action on $\D$ 
and it is transitive on $\C$. Hence the commutativity
condition in 
Proposition~\ref{p:equi} is automatic.
This raises the following question:
Let $\D$ be a thin labelled chamber complex of rank greater than 3
such that $\sigma_I \sigma_J = \sigma_J \sigma_I$ for all subsets $J
\subseteq I$.
Then is it true that $\D$ is a Coxeter complex ?
\end{remark}

\subsection{Connection with random walks} \label{ss:rw}
Bidigare, Hanlon, and
Rockmore \cite{\bhr} found a natural family of random walks
associated with hyperplane arrangements.  These walks were studied
further by Brown and Diaconis~\cite{\bd} and generalised to LRBs by
Brown~\cite{\ken}. One reason this development 
is exciting is that the walks admit a rather complete theory.   
We now describe the walk.

Using the projection operators, 
define the following walk on the chambers $\C$ of $\D$: 
If the walk is
in chamber $C$, choose a face $F$ of rank $j$ at random and move to the
projection $D=FC$.  (We assume here that the faces of rank $j$ are chosen
uniformly.)
Note that this random walk on chambers has a uniform stationary
distribution if and only if 
$\sigma_n \sigma_j = \sigma_j \sigma_n$.
Here $n$ is the rank of $\D$.

More generally, we can run a random walk on $\C$ using a probability
distribution $\{w_F\}_{F \in \F}$ on $\F$, the faces of $\D$.
Also if the complex is labelled then we can define the walk by
choosing instead a face $F$ of type $J$ at random.
Then the walk on chambers has a uniform stationary
distribution if and only if 
$\sigma_I \sigma_J = \sigma_J \sigma_I$.
However for the rest of the section, we will restrict to the unlabelled case.

\subsection{Some commutativity conditions}

Motivated by the above discussion, we look at a more general
commutativity problem. We work with a graded
LRB. All known examples of LRBs are graded. 
However we do not know whether this follows from the definition.
For now, we do not make any further assumptions.

Let $S$ be a \emph{LRB}, $L$ its associated lattice and $\supp\colon
S\onto L$ the support map.
For any $X \in L$, let
$S_{\le X} = \{ y \in S : \text{supp} \ y \le X \}.$
Then one can check that $S_{\le X}$ is also a LRB whose associated
lattice is the interval $[\zerohat,X]$ in~$L$.
Here $\zerohat$ is the support of the identity of $S$ and hence the
minimum element of $L$.
Also the set of chambers $\C$ consists precisely of
those elements of $S$ whose support is $\onehat$, 
the maximum element of $L$.

We say
that a LRB satisfies the commutativity condition for $i$ and $j$ if
\begin{equation} \tag{$C_{i,j}$}
\sigma_i \sigma_j = \sigma_j \sigma_i.
\end{equation}
Similarly we say that it satisfies the uniformity condition $(U)$ 
if for any i,j 
the number of times a chamber occurs in the product
$\sigma_i \sigma_j$ does not depend on the chamber we choose.
Loosely speaking, these conditions may be thought of as certain
symmetry conditions on the LRB.
We have already given some motivation for the condition $(C_{i,j})$. 
The motivation for condition $(U)$ becomes clear from the following
proposition. 

\begin{proposition}
Let $S$ be a LRB such that $S_{\le X}$ satisfies 
the uniformity condition $(U)$ for all $X \in L$. 
Then $S$ (and
hence all $S_{\le X}$) satisfies 
the commutativity condition $(C_{i,j})$ for all $i,j$. 
\end{proposition}
\begin{proof}
We first observe that 
$\supp(FG) = \supp(GF) = \supp F \vee \supp G$.
In particular, for $H=FG$, we have $\supp(F),\supp(G) \leq \supp(H)$.
Hence the coefficient of $H$ in $\sigma_i \sigma_j$ 
remains unchanged if we replace $S$ by the subsemigroup
$S_{\le \supp(H)}$. 
This shows that if $S$ satisfies $(C_{i,j})$ 
then so does $S_{\le X}$ for all $X \in L$.

Now we show that $S$ (and hence all $S_{\le X}$) 
satisfies $(C_{i,j})$ for all $i,j$.
We first check that for any $D \in \C$,
the coefficient of $D$ in $\sigma_i \sigma_j$ is same as 
the coefficient of $D$ in $\sigma_j \sigma_i$.
Since we assume that $S$ satisfies $(U)$, 
we simply need to check that 
$\abs{\{(F,G)\mid rk(F)=i, rk(G)=j, FG \ \text{is a chamber}\}} = 
\abs{\{(G,F)\mid rk(F)=i, rk(G)=j, GF \  \text{is a chamber}\}}$.
This is true since $FG$ is a chamber $\Leftrightarrow$ $GF$ is a chamber.
This is because of our earlier observation that 
$\supp(FG) = \supp(GF)$.

To do the general case, we need to check that for any $H \in \F$,
the coefficient of $H$ in $\sigma_i \sigma_j$ is same as 
the coefficient of $H$ in $\sigma_j \sigma_i$.
To do this computation, 
we replace $S$ by $S_{\le \supp(H)}$.
Note that $H$ is a chamber in $S_{\le \supp(H)}$ and by our assumption
$S_{\le \supp(H)}$ satisfies $(U)$. 
Hence we now apply the previous argument and the proposition is proved.
\end{proof}

\begin{corollary}
If $S$ is the Coxeter complex of type $A_{n-1}$ or $B_n$, then
$S$ satisfies $(C_{i,j})$ for all $i,j$.
\end{corollary}
\begin{proof}
The uniformity condition $(U)$ is automatic for a Coxeter complex $\Sigma$ 
because the Coxeter group acts transitively 
on the chambers of $\Sigma$. 
However to apply the above proposition, 
we require $(U)$ to hold for all $S_{\le X}$.
For types $A_{n-1}$ or $B_n$, 
given any $X \in L$, $S_{\le X}$ is again
isomorphic to a Coxeter complex of the same type 
but now with a smaller value of $n$.
Hence the Coxeter complex of type $A_{n-1}$ or $B_n$ satisfies $(U)$ for
all $S_{\le X}$ and we can apply the proposition.
\end{proof}

\begin{remark}
For types $A_{n-1}$ and $B_n$, the semigroups and the underlying
lattices can be made combinatorially explicit~\cite{\ken}.
The fact about $S_{\le X}$ being a Coxeter complex 
for types $A_{n-1}$ or $B_n$
that we used in the proof above is then 
quite elementary. For the case when $X$ is a hyperplane (wall), this
problem has been treated in general for any Coxeter complex by
Abramenko~\cite{\abra}.

A completely different proof of this corollary can be given using the
language of card shuffles~\cite{\riffle,\me}.
We also note that for $S$, the Coxeter complex of type $D_n$, the
subcomplexes $S_{\le X}$ are not even 
Coxeter complexes. This is the primary reason why the complex of
type $D_n$ fails the commutativity condition in general.
We have checked by explicit computation that 
the Coxeter complex of type $D_4$ does not satisfy
$(C_{i,j})$ for all $i,j$.
\end{remark}

\subsection{A low rank computation}
In contrast to the previous subsection, we get a less restrictive
result if we work with low ranks.

\begin{theorem}
Let $S$ be an oriented matroid of rank $3$. Then
$S$ satisfies $C_{1,2}$
$\Leftrightarrow$
$S$ satisfies $C_{1,3}$
$\Leftrightarrow$
$S$ satisfies $C_{2,3}$
$\Leftrightarrow$
$S$ is simplicial.
\end{theorem}
\begin{proof}
A part of this theorem, namely,
$S$ satisfies $C_{1,3}$
$\Leftrightarrow$
$S$ is simplicial
occurs as Proposition 1 in~\cite{\bbd}. 
It says that the random walk on $\C$ driven by a uniform probability
distribution on the vertices of $S$ has a uniform stationary
distribution iff $S$ is simplicial (see the discussion in
Section~\ref{ss:rw}).  
Here we will prove that
$S$ satisfies $C_{1,2}$
$\Leftrightarrow$
$S$ is simplicial.
The other cases are simpler and we leave them out.
The spirit of the computation is as in~\cite{\bbd}. 
Let $v$,$e$ and $f$ be the number of vertices, edges and faces
respectively in $S$.

We need to look at $\sigma_1 \sigma_2 = \sigma_2 \sigma_1$.
Note that a vertex can never occur in either product.
And the number of times an edge appears on either side is the same.
(This is because we can restrict to the support of that edge, which is
just a pseudoline.)
So the non-trivial case is that of chambers. We fix a chamber $D$.
Denote the chamber exactly opposite to $D$ by $\overline{D}$.
We count the number of times that either $D$ or $\overline{D}$ occurs in
$\sigma_1 \sigma_2$ and $\sigma_2 \sigma_1$.
(It will be evident from the computation as to why we group $D$ and
$\overline{D}$ together. Note that due to the antipodal symmetry of the
arrangement, we do not lose anything essential by doing this.)

\begin{figure}[hbt]
\centering
\begin{tabular}{c@{\qquad}c}
\mbox{\epsfig{file=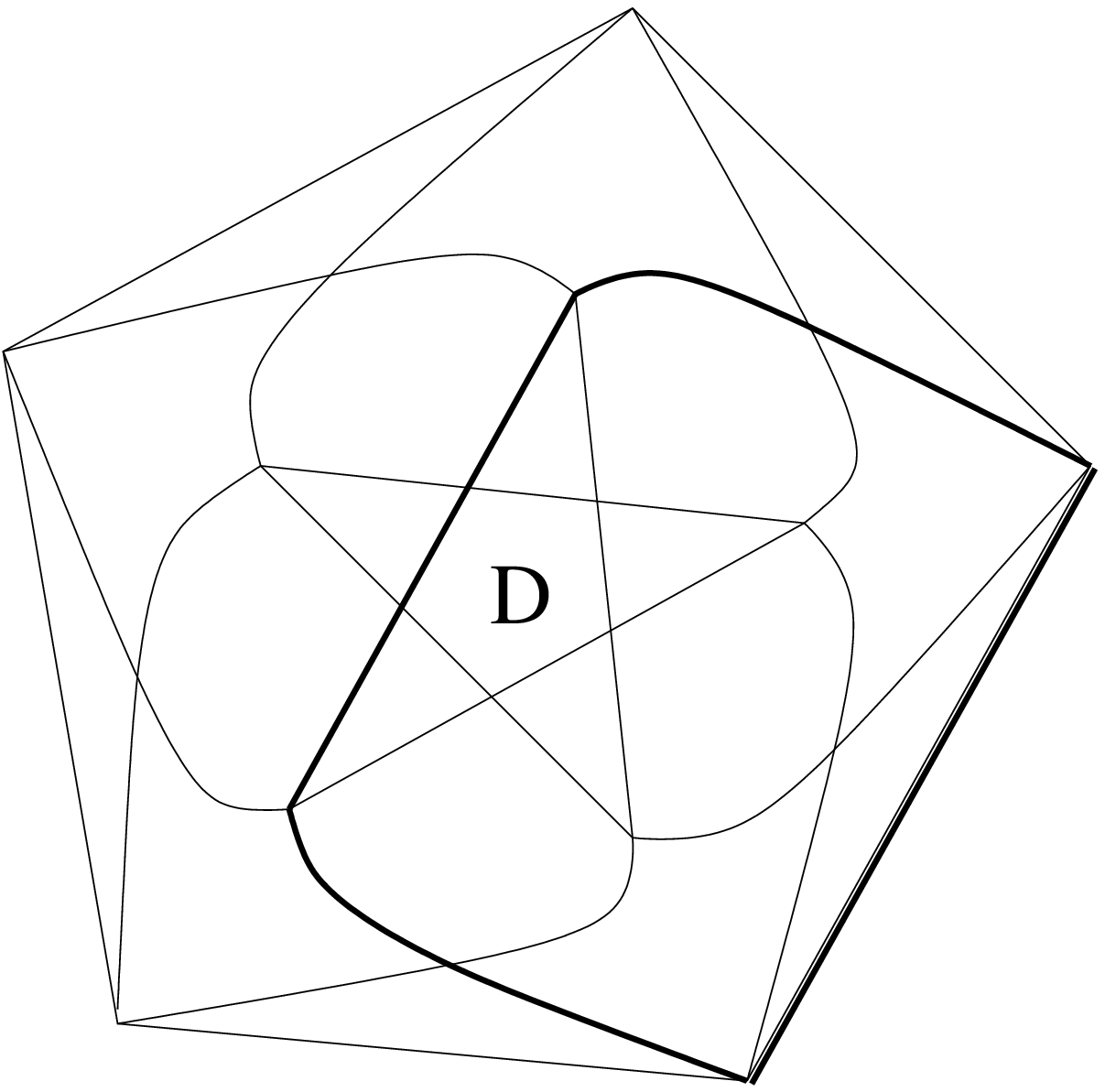,height=3cm,width = 3cm}}
\end{tabular}
\caption{The walls of a $5$-gon.  }
\label{pentagon}
\end{figure}

If $D$ is a $k$-gon then consider the $k$ walls of $D$
(and hence $\overline{D}$). 
We call an edge of $S$ an \emph{interior} edge if it does not lie on any of
these $k$ walls.
Figure~\ref{pentagon} illustrates the case
when $k=5$. One of the walls has been highlighted. 
The only chamber that is not seen in the picture is the pentagon
$\overline{D}$. It is on the backside bounded by the five
segments that bound the figure.

We first look at $\sigma_1 \sigma_2$ and count the number of times
that either $D$ or $\overline{D}$ occurs in this product; that is,
$\abs{\{(F,G)\mid FG=D \  \text{or} \  FG=\overline{D}, rk(F)=1, rk(G)=2
\}}$.
We use the geometry of the arrangement as illustrated in the figure.
Note that $F$ must be a vertex of $D$ or $\overline{D}$. 
We consider two cases.
In the first case, let $G$ be an interior edge or an edge of either
$D$ or $\overline{D}$.
Then there are $k-2$ choices for $F$. 
The individual contributions of $D$ or $\overline{D}$ depend on the
location of the edge $G$, however the net contribution is always $k-2$.
In the second case, $G$ lies on one of the $k$ walls but is not a face 
of either $D$ or $\overline{D}$.
Now there are $k-3$ choices for $F$.
Hence the required count is $(k-2)e - (\text{number of edges on all the
walls}) + 2k$.

Now we look at $\sigma_2 \sigma_1$. We want to count 
$\abs{\{(G,F)\mid GF=D \  \text{or} \  GF=\overline{D}, rk(F)=1, rk(G)=2
\}}$. 
The only $G$'s that we need consider are the edges of $D$ and $\overline{D}$.
Given a wall $H$ of $D$ (and $\overline{D}$), the edge of $D$
that lies on $H$ and its opposite edge (which is an edge of $\overline{D}$)
together can pair with all vertices except the ones that lie on $H$.
Summing up over all the $k$ walls, we get the count
to be 
$kv - (\text{the sum of the no. of vertices on each of the $k$ walls})$.

Now we compare the counts for $\sigma_1 \sigma_2$ and $\sigma_2 \sigma_1$.
Since the no. of vertices is same as the no. of edges on each
wall (since it is a pseudoline), we get
$\sigma_1 \sigma_2 = \sigma_2 \sigma_1
\Leftrightarrow
(k-2)e + 2k = kv
\Leftrightarrow
kf =2e
$.
The last equality holds if and only if the arrangement is made entirely of
$k$-gons.
However, by Levi's
theorem~\cite[p.~25]{\gruna}),
any arrangement must have at least one triangle. 
Hence we get $k=3$.
\end{proof}

\section{The opposite or duality axiom}

Let $\D$ be a simplicial complex of rank $n$.
Under suitable hypotheses on $\D$ 
(for example, if $\D$ triangulates a sphere),
the $h$ vector of $\D$ satisfies the linear relations $h_i=h_{n-i}$ for
$0 \leq i \leq n$.
Similarly if $\D$ is a labelled simplicial complex with label set $I$
then under suitable hypotheses the flag $h$ vector of $\D$ satisfies
the linear relations 
$h_J = h_{I \setminus J}$,
for $J \subseteq I$.
The former are called the \emph{Dehn-Sommerville} equations 
and the later the
\emph{generalised Dehn-Sommerville} equations. 
From now on, we will just write them as DS for short.
For more information see \cite[Section 3.14]{\stanley} and the
references therein.

The DS equations are a form of Alexander duality.
Our goal for this section is to formulate an axiom
relevant to this notion of duality.
As motivation we begin with a result that
introduces the concept of shelling reversal.

\vanish{
Let $\hat{P}$ be a ranked poset of rank $n+1$ with $0$ (resp. $1$) as
the minimal (resp. maximal) element.
The flag (order) complex of $P=\hat{P}-\{0,1\}$, which we call $\D(P)$, is a labelled
simplicial complex of rank $n$. 
The flag $h$ vector of $\hat{P}$ is by definition the flag $h$ vector
of the flag (order) complex of $\D(P)$.
Under suitable hypotheses on $\hat{P}$ , all linear relations among
the entries of its flag $h$ vector are given by the generalised
Dehn-Sommerville (DS for short) equations,
$h_J = h_{I \setminus J}$,
where $I$ is the full set of active ranks $\{1,2,\ldots,n\}$ and $J \subseteq I$.
}

\begin{theorem} \label{t:dual}
Let $\D$ be a shellable complex such that every facet lies in at least
2 chambers. Then
$\D$ is thin 
$\Rightarrow$
$\abs{\{D \mid R(D) = D\}}=1$
$\Leftrightarrow$
$\D$ is homotopy equivalent to a sphere.

If $\D$ is thin then it satisfies the DS equations.
\end{theorem}
\begin{proof}
The first part is proved in~\cite[Theorem 1.5]{\bjorner}.
We prove the second part. 
Also we assume that $\D$ is labelled.
If it is not then the same proof gives us the DS equations instead of
the generalised DS equations.
The proof is based on the idea of shelling reversal.
 
Let $\leq_S$ be a shelling order of $\D$, and let 
$\geq_S$ be the reverse order defined by
$D \geq_S E \Leftrightarrow D \leq_S E$;
that is, the linear order on $\C$ 
specified by $\geq_S$ is exactly reverse to that of $\leq_S$.
Assume for the moment that $\geq_S$ also gives a shelling of $\D$.
Now if we
let $R$ (resp. \reflectbox{R}) be the restriction map of 
$\leq_S$ (resp. $\geq_S$)
then it follows that 
$$\text{type of} \ R(D) = I \setminus \ \text{type of} \ \reflectbox{R}(D).$$
This along with the fact that 
$h_J = \abs{\{D \mid R(D) \  \text{has type} \  J\}}$
gives us the generalised DS equations 
and proves the second part.

\begin{figure}[hbt]
\centering
\begin{tabular}{c@{\qquad}c}
\mbox{\epsfig{file=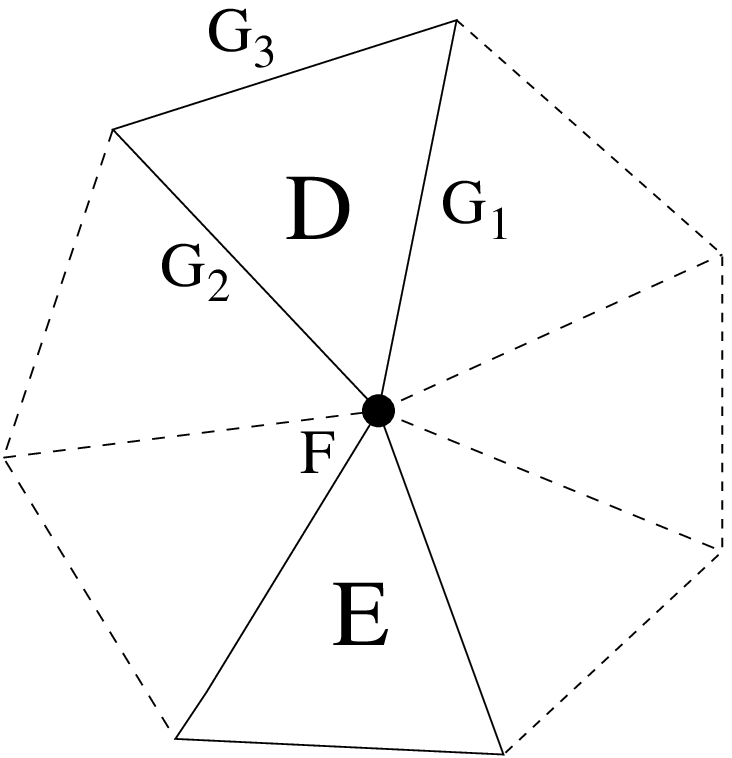,height=3cm,width = 3cm}}
\end{tabular}
\caption{An example with $\F_D=\{G_1,G_2\}$ and $\F^{\prime}_D=\{G_3\}$.}
\label{opposite}
\end{figure}

\noindent
Hence we only need to show that $\geq_S$ gives a shelling of $\D$.
We assumed that $\leq_S$ is a shelling of $\D$.
This says that
$D \cap (\cup_{E <_S D} E) = \cup_{G \in \F_D} G$,
where $\F_D$ is a subset of the set of facets of $D$.
Now let $\F^{\prime}_D$ be the set of those facets of $D$ 
that are not in $\F_D$. 
We will be done once we show the following.

\medskip
\noindent
Claim: $D \cap (\cup_{E >_S D} E) = \cup_{G \in \F^{\prime}_D} G$.\\
Proof of claim. $(\supseteq)$ This is clear.\\
$(\subseteq)$
Let $F$ be a face of $D$ such that $F \nsubseteq RHS$.
We show that $F \nsubseteq LHS$.
Let $E \geq F$ such that $E \not= D$.
We need to show that $E <_S D$.
To put in words, we need to show that among all the chambers in
$\D_{\geq F}$, $D$ is the last chamber to be shelled
by the order $\leq_S$.
By Lemma~\ref{l:sh}, we know that $\D_{\geq F}$, or equivalently
$lk(F,\D)$, is shellable with the induced order from $\leq_S$.
The chambers of $lk(F,\D)$ are in 1--1 correspondence with the chambers
of $\D$ that contain $F$.
Since we assume that $\D$ is thin, it follows that $lk(F,\D)$ is
also thin. 
Let $R_{F}$ be the restriction map of the induced shelling on $lk(F,\D)$.
Let $D_F$ be the chamber of $lk(F,\D)$ that corresponds to $D$.
We need to show that $D_F$ is the last chamber to be shelled in
$lk(F,\D)$. Since $lk(F,\D)$ is thin, applying the first part of the
theorem, we need to show that $R_{F}(D_F)=D_F$; also see comment after
Proposition~\ref{p:shell}.
However this follows by our assumption that $F \nsubseteq RHS$, which says that every facet of
$D$ that contains $F$ belongs to $\F_D$. 
This proves the claim.
\end{proof}

\begin{remark}
For a sketch of the same argument in the context of polytopes, 
see~\cite[pg 252]{\zie}.
Also for the first part of the theorem, we do not know whether
thinness is equivalent to the other two conditions.
\end{remark}

\subsection*{The opposite axiom}
Let $\D$ be a thin shellable complex.
Let $C$ and $D$ be the first and last elements of a shelling $\leq_S$
of $\D$.
Then $C$ and $D$ are the unique elements satisfying $R(C)=\emptyset$ and
$R(D)=D$, where $R$ is the restriction map of $\leq_S$.
As we saw in the proof of Theorem~\ref{t:dual}, the shelling $\leq_S$
is reversible and for the reverse shelling $\geq_S$, the roles of $C$
and $D$ get interchanged. Hence we would like to think of $C$ and $D$
as chambers opposite to each other.

Motivated by this discussion, we add an \emph{opposite} axiom to the
axiomatic setup of Section~\ref{s:axiom}.
Of course, the main case of interest is when $\D$ is thin.
We give three equivalent definitions 
which correspond to the projection, restriction and shelling cases
respectively.
Their equivalence will be proved in Theorem~\ref{t:opposite}.

\begin{figure}[hbt]
\centering
\begin{tabular}{c@{\qquad}c}
\mbox{\epsfig{file=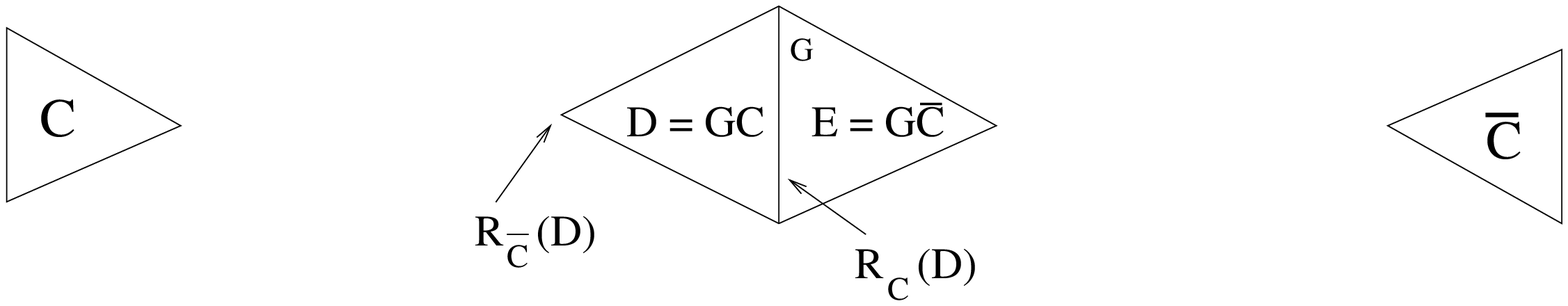,height=2.4cm,width = 11cm}}
\end{tabular}
\caption{How the opposition map works.}
\label{opp}
\end{figure}

\noindent
Let $^{-}:\C \rightarrow \C$ be a map such that
\begin{itemize}
\item[(P4)]
For any facet $G$ of $\D$ and any $C \in \C$, we have
$GC \not= G\overline{C}$.
 
\item[(R4)] 
For all $C, D \in \C$, we have 
type of $R_C(D) = I \ \setminus$ type of $R_{\overline{C}}(D)$.

\item[(S4)]
For all $C \in \C$, the partial orders $\leq_C$ and $\leq_{\overline{C}}$
are dual to each other.
\end{itemize}
If $\D$ is not labelled then axiom $(R4)$ is to be interpreted as
saying that $R_C(D)$ and $R_{\overline{C}}(D)$ are complementary faces
of $D$. 
If $\D$ is thin then axiom $(P4)$ says that for a facet $G$, $GC$ and
$G\overline{C}$ are the two distinct chambers containing $G$.
Each of these axions implies that 
$^{-}:\C \rightarrow \C$ 
is an involution of $\D$. 
Also, using these axioms, ${\overline{C}}$, ``the
chamber opposite to $C$'', can be characterised as
follows in the three cases. 
\begin{itemize}
\item[P.]
$\overline{C}$ is the unique chamber such that 
$FC \not= \overline{C}$ for any $F < \overline{C}$.
 
\item[R.] 
$\overline{C}$ is the unique chamber satisfying 
$R_C(\overline{C}) = \overline{C}$.

\item[S.]
$\overline{C}$ is the unique maximal element in the partial order $\leq_C$.
\end{itemize}

\begin{remark}
Recall that 
$h_J(D) = \abs{\{C \mid R_C(D) \  \text{has type} \  J\}}$. 
We claim that if $\D$ satisfies the opposite axiom then
$h_J(D) = h_{I \setminus J}(D)$. 
This is because by axiom $(R4)$, if $C$ contributes to 
$h_J(D)$ then its opposite $\overline{C}$ contributes to 
$h_{I \setminus J}(D)$. 
Now averaging over all $D \in \C$ leads to the generalised DS equations.
(Of course, we already knew from Theorem~\ref{t:dual} that $\D$
satisfied the generalised DS equations.)
Hence we may regard $h_J(D) = h_{I \setminus J}(D)$ as a finer version
of these equations.
In the unlabelled case, the same reasoning gives a finer version of the DS equations
instead of the generalised DS equations.
\end{remark}

\begin{theorem} \label{t:opposite}
Let $\D$ be a thin chamber complex that satisfies our (previous)
axioms.
Then
$\Delta$ satisfies $(P4)$ $\Leftrightarrow$ $\Delta$ satisfies $(R4)$ $\Leftrightarrow$ $\Delta$ satisfies $(S4)$.
\end{theorem}
\begin{proof}
The complex $\D$ is shellable since it satisfies our earlier axioms.
And a shellable complex is automatically gallery connected. 
So really the only additional assumption on $\D$ is that of thinness.
To go from one set of axioms to another, we again use the formal
connections outlined in Section~\ref{subs:ca}.

\subsection*{$(P4) \Rightarrow (R4)$}
Let $D \in \C$. Apply axiom $(P4)$ to every facet of $D$. 
By the definition of $R_C$ and $R_{\overline{C}}$, we get that the
type of $R_C(D) = I \ \setminus$ type of $R_{\overline{C}}(D)$.

\subsection*{$(R4) \Rightarrow (S4)$}
For a thin shellable complex,
Lemma~\ref{l:comp} gives a simpler description of $\leq_C$.
It is the
transitive closure of the relation:
$E \leq_C D$ if $E$ is adjacent to $D$ and $R_C(E) \leq D$.

Let $D$ and $E$ be adjacent chambers and $D\not=E$.
Let $G$ be the common facet.
By thinness of $\D$ and axiom $(R2)$, 
either $R_{\overline{C}}(D) \leq G$ or
$R_{\overline{C}}(E) \leq G$.
This implies that either $R_{\overline{C}}(D) \leq E$ or
$R_{\overline{C}}(E) \leq D$.

In order to prove axiom $(S4)$, 
we need to show that 
$D \leq_C E \Leftrightarrow E \leq_{\overline{C}} D$.
We show that
$D \leq_C E \Rightarrow E \leq_{\overline{C}} D$.
The other implication follows by symmetry.
Suppose that $D$ and $E$ are adjacent and $R_C(D) \leq E$.
Then by axiom $(R4)$,
we have $R_{\overline{C}}(D) \not \leq E$.
Hence $R_{\overline{C}}(E) \leq D$
which implies $E \leq_{\overline{C}} D$.

\subsection*{$(S4) \Rightarrow (P4)$}
Let $G$ be any facet of $\D$. Since $\D$ is thin, there are exactly
two chambers in $\D_{\geq G}$.
Axiom $(S1)$ says that for any $C \in \C$, the partial order $\leq_C$
restricted to $\C_{\geq G}$ has a minimum.
So the two chambers in $\D_{\geq G}$ get ordered by $\leq_C$.
Now axiom $(S4)$ says that the order obtained using
$\leq_{\overline{C}}$ is reverse of the one got using $\leq_C$ and
hence axiom $(P4)$ follows by the definition
of projection maps.
\end{proof}

\noindent
A simplicial hyperplane arrangement clearly satisfies 
axiom $(P4)$.
Axiom $(R4)$ can also be checked directly
using the discussion in Section~\ref{ss:example}.
We record this formally as a simple corollary.

\begin{corollary}
A simplicial hyperplane arrangement satisfies the opposite
axioms $(P4)$, $(R4)$ and $(S4)$.
\end{corollary}

\section{Duality in Buildings} \label{s:dbuild}

In the previous section, we axiomatized duality in a way that was
relevant to thin chamber complexes.
Now we prove a somewhat isolated result on duality in buildings
and illustrate it with the building of type $A_{n-1}$.
A suitable generalisation of the opposite axiom to buildings may
be possible.
An appropriate setting may be gated chamber complexes
of spherical type~\cite{\muhlherr}.
Also the duality result on buildings (Theorem~\ref{t:duality}) 
may generalise to Moufang 
complexes~\cite[Section 5]{\schmid}.

We first set up some notation.
Let $\D$ be a finite (spherical) building and $\A$ be the complete set of
apartments of $\D$.
For chambers $C,D$, we denote by $\A_C$ (resp. $\A_{C,D}$) the set of
apartments in $\A$ that contain $C$ (resp. $C$ and $D$).
Also, for a chamber $C$, let $C^{op}$ denote the set of all chambers
that are opposite to $C$. 
Observe that there is a bijection between the sets $C^{op}$ and
$\A_C$ as follows.
We associate to $\overline{C} \in C^{op}$, the apartment $\Sigma \in
\A_C$, which is the convex hull of $C$ and $\overline{C}$.
Now let $\A^{\prime}_{C,D} = \{E \in C^{op} \mid D \ \text{lies in
the apartment determined by} \ C \ \text{and} \ E \}$.
Observe that the bijection between $\A_C$ and
$C^{op}$ restricts to a bijection between 
the sets $\A_{C,D}$ and $\A^{\prime}_{C,D}$.
Also we let $\rho_{\Sigma,C} : \D \rightarrow \Sigma$ denote the usual
retraction. 

\begin{lemma} \label{l:duality}
Let $C$ and $\overline{C}$ be a pair of opposite chambers in $\D$ 
and let $D$ be any chamber in the apartment determined by 
$C$ and $\overline{C}$. Then
$\abs{\A_{C,D}} \abs{\A_{D,\overline{C}}} = \abs{\A_{D}}.$
\end{lemma}
\begin{proof}

\begin{figure}[hbt]
\centering
\begin{tabular}{c@{\qquad}c}
\mbox{\epsfig{file=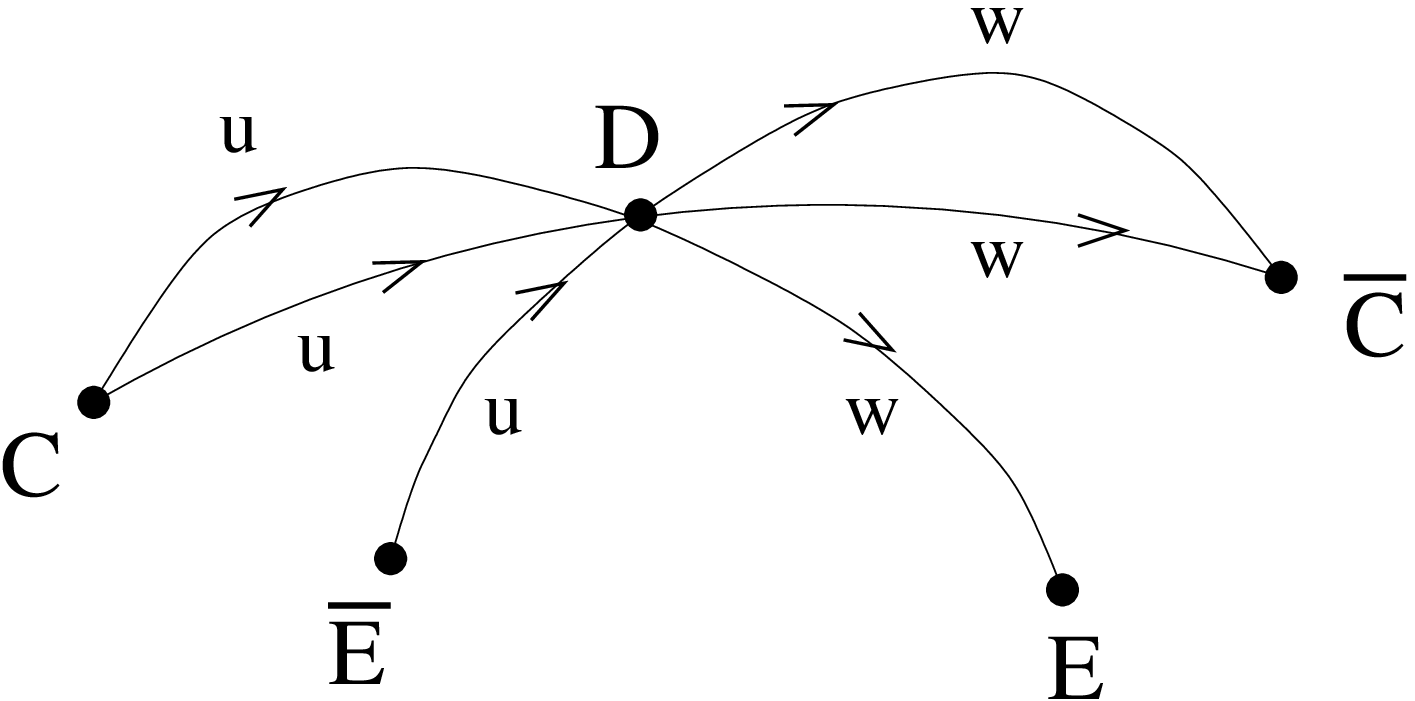,height=2.5cm,width = 6cm}}
\end{tabular}
\caption{}
\label{building}
\end{figure}

Let $\delta : \C \times \C \rightarrow W$ be the $W$-valued distance
function on $\C$, where $W$ is the Coxeter group of $\D$.
Since $\D$ is spherical the Coxeter group $W$ is finite.
Denote the longest word in $W$ by $w_0$.
Let $C,D$ and $\overline{C}$ be as in the statement of the lemma.
Also let $\delta(C,D)=u$ and $\delta(D,\overline{C})=w$. Then $uw=w_0$
since $C$ and $\overline{C}$ are opposite 
and $C,D$ and $\overline{C}$ lie in the same apartment.
Note that $\A^{\prime}_{C,D} = \{ E \mid \delta(D,E) = w \}$.
This is because $l(uw)=l(u)+l(w)$.
We now define a bijection $\A^{\prime}_{C,D} \times \A^{\prime}_{D,\overline{C}} \rightarrow \A_{D}$
as follows.
If $E \in \A^{\prime}_{C,D}$ and $\overline{E} \in \A^{\prime}_{D,\overline{C}}$
then $\delta(\overline{E},E)=uw=w_0$ again because $l(uw)=l(u)+l(w)$.
Hence $E$ and $\overline{E}$ are opposite.
So we map $(E,\overline{E})$ to the unique apartment containing $E$
and $\overline{E}$. 
The map is clearly a bijection.

\end{proof}

Let $\D$ be a finite \emph{irreducible} Moufang building. 
By irreducible, we mean that the Coxeter group $W$ associated to $\D$
is irreducible. We will not state the Moufang condition here; 
see~\cite[Chapter 6]{\ronan}. 
We only mention that all buildings of rank greater or equal to 3 are
Moufang. Hence it is a restriction only when $rk(\D)=2$.
By a theorem of Tits~\cite[Theorem 11.4]{\tits}, 
if $\D$ is irreducible and Moufang then it is the building of an
algebraic group $G$ over a finite field $\fq$.
The group $G(\fq)$ acts by simplicial type-preserving automorphisms on
$\D$. Also the action is transitive on pairs $(\Sigma,C)$ consisting of an
apartment $\Sigma \in \A$ and a chamber $C \in \Sigma$.
In particular, the stabiliser of a fixed chamber $C$, which we denote
$\Stab(C)$ acts transitively on $\A_C$.
It follows from the properties of buildings that for $g \in \Stab(C)$,
$\rho_{\Sigma,C}(gD) = \rho_{\Sigma,C}(D)$ for any chamber $D \in \C$.

It is known that for any chamber $C$ of $\D$, the number of chambers
adjacent to $C$ along a fixed facet is a power of $q$.
Here $q$ is the cardinality of the finite field $\fq$.
It follows that for chambers $C,D$, the number $\abs{\A_{D}}$ (and also
$\abs{\A_{C,D}}$) is a power of $q$. 
Using this fact, $h_J(\D)$, the flag $h$ vector of $\D$, defined in
Section~\ref{s:tss}
can be written
as a polynomial in $q$. 
In the
course of proving Theorem~\ref{t:duality},
we will recall the definition of $h_J(\D)$ and give 
a precise definition of this polynomial.
To emphasize the dependence on $q$, we write $h_J(q)$ instead of $h_J(\D)$.
It will also follow that
if $\Sigma$ is any apartment of $\D$ then
$h_J(\Sigma) = h_J(q) \mid_{q=1}$.
The geometric significance of $h_J(\D)$ was explained in
Propositions~\ref{p:type} and~\ref{p:beta}.
It counts the number of spheres in the homotopy type of $\D_J$,
the type selected subcomplex of $\D$.
Now we prove the main result of this section.
It gives the precise relation between the homotopy types of
$\D_J$ and $\D_{I \setminus J}$.

\begin{theorem} \label{t:duality}
Let $\D$ be a finite, irreducible and Moufang building.
Then with the notation as above $h_J(q) = h_I(q) h_{I \setminus J}(q^{-1})$.
\end{theorem}
\begin{proof}
We start with a philosophical comment.
``As $q \rightarrow 1$, the building $\D$ 
degenerates to a single apartment $\Sigma$.
The equation of the theorem then reduces to the generalised DS equation
$h_J(\Sigma) = h_{I \setminus J}(\Sigma)$.''
Hence to prove the theorem,
we will go in the opposite direction.
We will start with the relation
$h_J(\Sigma) = h_{I \setminus J}(\Sigma)$. Then
for a fixed pair of opposite chambers $C, \overline{C} \in \Sigma$,
we will use the retractions
$\rho_{\Sigma,C}$ and $\rho_{\Sigma,\overline{C}}$
to generalise this relation to $\D$.
We now start with the formal proof.

We fix an apartment $\Sigma$ and a chamber $C$ that lies in it.
Recall that $R_C(D)$ is the smallest face $F$ of $D$ such that $FC=D$
(Section~\ref{subs:ca}).
Now set $H_J(\D) = \{D \in \D \mid R_C(D) \  \text{has type} \  J\}$. 
Hence by definition, the cardinality of this set is $h_J(\D)$.
Similarly let $H_J(\Sigma,C) = \{D \in \Sigma \mid R_C(D) \
\text{has type} \  J\}$.  
We make a few observations.
The restriction map $R_C$ for $\D$ when restricted to the apartment
$\Sigma$ gives the restriction map for $\Sigma$.
This follows, for example, from the fact that the projection maps
for $\D$ are defined using the projection maps for the apartments.
Hence the cardinality of $H_J(\Sigma,C)$ is $h_J(\Sigma)$.

Next note that 
$\Sigma$ is a Coxeter complex and 
hence corresponds to a simplicial hyperplane arrangement.
So it satisfies the opposite axiom $(R4)$;
see the Corollary to Theorem~\ref{t:opposite} and the preceding comments.
Hence applying axiom $(R4)$ to it, we get $H_J(\Sigma,C) = H_{I \setminus
J}(\Sigma,\overline{C})$.
This is a more precise set theoretic version 
of the generalised DS equation.
The cardinalities of the two sets are $h_J(\Sigma)$ and $h_{I
\setminus J}(\Sigma)$ respectively and do not depend on $C$ and $\overline{C}$.
We next claim that 
\begin{equation*}
H_J(\D) = \sqcup_{D \in H_J(\Sigma,C)} \rho^{-1}_{\Sigma,C}(D).
\end{equation*}
To see this, note that
$\rho_{\Sigma,C}(F) C = \rho_{\Sigma,C}(F C)$ for any face $F$ of $\D$.
Hence $F=R_C(D) \Leftrightarrow \rho_{\Sigma,C}(F) =
R_C(\rho_{\Sigma,C}(D))$ 
which proves the claim.
Let $\overline{C}$ be the chamber opposite to $C$ in $\Sigma$.
Then $\rho^{-1}_{\Sigma,C}(\overline{C}) = C^{op}$, 
the set of chambers opposite to $C$.
Also note that
$h_I(\Sigma) = \abs{H_I(\Sigma,C)} = 1$
with $\overline{C}$ as the only element of $H_I(\Sigma,C)$.
Applying the claim for $J=I$,
we get
$h_I(\D) =$ the number of spheres in the homotopy
type of $\D = \abs{C^{op}} = \abs{\A_C}$.
We now have the following identities.

\begin{equation*}
h_J(q) = \sum_{D \in H_J(\Sigma,C)} \abs{\rho^{-1}_{\Sigma,C}(D)} =
\sum_{D \in H_J(\Sigma,C)} \frac{\abs{\A_D}}{\abs{\A_{C,D}}}.
\end{equation*}

\begin{equation*}
h_{I \setminus J}(q) = \sum_{D \in H_{I \setminus J}(\Sigma,\overline{C})} \abs{\rho^{-1}_{\Sigma,\overline{C}}(D)} =
\sum_{D \in H_J(\Sigma,C)} \frac{\abs{\A_D}}{\abs{\A_{D,\overline{C}}}}.
\end{equation*}
The first equality is clear.
For the second equality we apply Tits theorem 
to get a group $G(\fq)$ associated to $\D$.
The group serves two purposes.
It gives us the number $q$ and it acts nicely on $\D$.
As already noted $\Stab(C)$ (resp. $\Stab(\overline{C})$) 
acts transitively on $\A_C$ (resp. $\A_{\overline{C}}$) 
and is compatible with the 
retraction $\rho_{\Sigma,C}$ (resp. $\rho_{\Sigma,\overline{C}}$). 
This along with the fact $\abs{\A_C}=\abs{\A_D}=\abs{\A_{\overline{C}}}$
gives us the second equality.
The formulas show that $h_J(q)$ can be written as a polynomial in $q$ as
claimed earlier.
Furthermore due to the group action, it does not depend on the choice of
$\Sigma$ and $C$. It is also clear that 
if $\Sigma$ is any apartment of $\D$ then
$h_J(\Sigma) = h_J(q) \mid_{q=1}$.

Recall that $h_I(q) = \ \text{the number of spheres in the homotopy
type of} \ \D = \abs{\A_D}$.
Now if we replace $q$ by $q^{-1}$ in the second equation and multiply by
$h_I(q)$, we get

$$h_I(q) h_{I \setminus J}(q^{-1}) = \abs{\A_D} 
\sum_{D \in H_J(\Sigma,C)} \frac
{\abs{\A_{D,\overline{C}}}}{\abs{\A_D}} = 
\sum_{D \in H_J(\Sigma,C)} {\abs{\A_{D,\overline{C}}}}.$$
The theorem now follows from Lemma~\ref{l:duality}.
\end{proof}

\begin{example} \label{e:dbuild}
Building of type $A_{n-1}$:
Let $V$ be the $n$-dimensional vector space $\fq^n$, where $\fq$ is
the field with $q$ elements.  
Let $\mathcal{L}_{n}$ be the lattice of subspaces
of $V$ under inclusion. 
Also let $\mathcal{B}_{n}$ be the Boolean lattice.
Note that the flag (order) complex $\D(\mathcal{B}_{n})$ 
is the Coxeter complex of type
$A_{n-1}$ and corresponds to the braid arrangement. Also note that a
choice of a basis for $V$ gives an 
embedding of $\mathcal{B}_{n}$ into $\mathcal{L}_{n}$.
The building of type $A_{n-1}$ over $\fq$ 
is defined simply as the flag (order) complex
$\D(\mathcal{L}_{n})$. The subcomplexes $\D(\mathcal{B}_{n})$, for various embeddings $\mathcal{B}_{n} \into \mathcal{L}_{n}$,
play the role of apartments.
The algebraic group $G(\fq)$ in this case can be taken to be $GL_n(\fq)$ or 
$SL_n(\fq)$.

For any chamber $D$, one may directly check that 
$\abs{\A_{D}} = h_I(q) = q^{\binom{n}{2}}$.
Now fix an apartment $\Sigma$ corresponding to a basis
$e_1,e_2,\ldots,e_n$ of $V$. 
Let $C = [e_{1}] < [e_{1},e_{2}] < \ldots < [e_{1},\ldots,e_{n}]$.
Then the chamber opposite to $C$ in $\Sigma$, namely 
$\overline{C} = [e_{n}] < [e_{n},e_{n-1}] < \ldots < [e_{n},\ldots,e_{1}]$.
If we let $D = [e_{i_1}] < [e_{i_1},e_{i_2}] < [e_{i_1},\ldots,e_{i_n}]$ be
any chamber in $\Sigma$ then a simple counting argument shows that
$$\abs{\A_{C,D}} = q^{\abs{\{(j,k) \mid j<k \ \text{and} \ i_j < i_k\}}}$$
$$\abs{\A_{D,\overline{C}}} = q^{\abs{\{(j,k) \mid j<k \ \text{and} \ i_j > i_k\}}}.$$
Hence we can see directly that 
$\abs{\A_{C,D}} \abs{\A_{D,\overline{C}}} = \abs{\A_{D}}.$
This proves Lemma~\ref{l:duality} for this example.
To unwind Theorem~\ref{t:duality}, we need to understand
how retractions work. We leave that out and instead just give some
explicit computations that illustrate the theorem. 
For type $A_3$, 
$h_{\emptyset}=1$, $h_{\{1\}}=h_{\{3\}}=q(1+q+q^2)$,
$h_{\{2\}}=q(1+2q+q^2+q^3)$,
$h_{\{1,2\}}=h_{\{2,3\}}=q^3(1+q+q^2)$,
$h_{\{1,3\}}=q^2(1+q+2q^2+q^3)$,
$h_{\{1,2,3\}}=q^6$.
We mention that these polynomials can also be described 
using descent sets without any
reference to retractions, see \cite[Theorem 3.12.3]{\stanley}.
\end{example}

\section{Future prospects}

In this paper we presented a theory of projection maps and compatible
shellings. We now suggest some problems for future study.

\subsection*{The space of all shellings}
As we saw, compatible shellings on a complex $\D$ consist
of compatible partial orders $\leq_C$, one for every chamber $C$ of
$\D$. 
Using this information, we may construct a geometric object $S(\D)$
containing a lot of shelling information about $\D$. It would be
worthwhile to make this idea precise since the object $S(\D)$ would be
important for studying the space of all shellings of $\D$. In the
example of Section~\ref{e:nmetric}, a candidate for $S(\D)$ is a line
with vertices at the integer points.  The group $\Z$ acts on $S(\D)$
by translations. The complex $\D$ is a quotient of $S(\D)$ by the
subgroup $n\Z$ of $\Z$ whose generator shifts a vertex $n$ units to the
right.

In this regard, the following formalism could be useful.
Let $\D$ be a pure simplicial complex that satisfies our axioms.
Let $\T$ be the space of all partial orders on $\C$ 
that extend $\leq_C$ for some $C \in \C$.
Then $\T$ is a poset with the order relation given by extension
with the partial orders $\leq_C$ as the minimal elements.
We also think of $\C$ as a poset with the trivial partial order.
Then $\T$ and $\C$ are posets related by a Galois connection.
Denote a typical element of $\T$ by $\leq_S$.
We are slightly abusing notation
because in our earlier usage, 
$\leq_S$ was always a linear order.
The map $\T \rightarrow \C$
sends $\leq_S$ to the chamber that appears first in $\leq_S$.
And the map $\C \rightarrow \T$
sends $C$ to $\leq_C$.

\subsection*{Other questions}
In contrast to the situation in Section~\ref{s:metric},
a shelling of a simplicial complex $\D$
may have nothing to do with its gallery metric;
see the example of Section~\ref{e:nmetric}.
We ask for other non-metrical examples of $\D$
that satisfy our axioms.

While we explained some tiny applications and connections of our
theory to the flag $h$ vector, random walks and buildings, other
applications of this circle of ideas still remain to be seen.  As a
possibility, we suggest that a general framework to study Solomon's
descent algebra~\cite{\solomon} would be to start with a labelled simplicial
complex $\D$ that satisfies our axioms. If $\D$ is a Coxeter complex
then we would recover the usual descent algebra;
see~\cite{\bid,\ken,\me} for more information on this geometric way of
thinking.

As shown by the examples in Section~\ref{s:proj}
our theory is closely related to LRBs.
We ask if this connection can be clarified.
For example, we may ask the converse question to 
the discussion in Section~\ref{s:LRB}.
Given a simplicial complex that satisfies our axioms,
how close is it to being a simplicial LRB 
(see the remark in Section~\ref{e:nmetric})?

\subsection*{More examples}
More examples of the theory that we presented may be
possible.  A possibility already mentioned is that of shellable
complexes with group actions. Another possibility is that of
simplicial polytopes equipped with a suitable class of shellings like
line shellings. However a better alternative seems to be to generalise
the theory itself. We conclude by suggesting some approaches through
examples.

$(i)$
Recall that in the metric setup, $FC$ was the chamber containing $F$
closest to $C$.
What happens if there is no closest chamber?
We present a building-like example where this occurs.

\begin{figure}[hbt]
\centering
\begin{tabular}{c@{\qquad}c}
\mbox{\epsfig{file=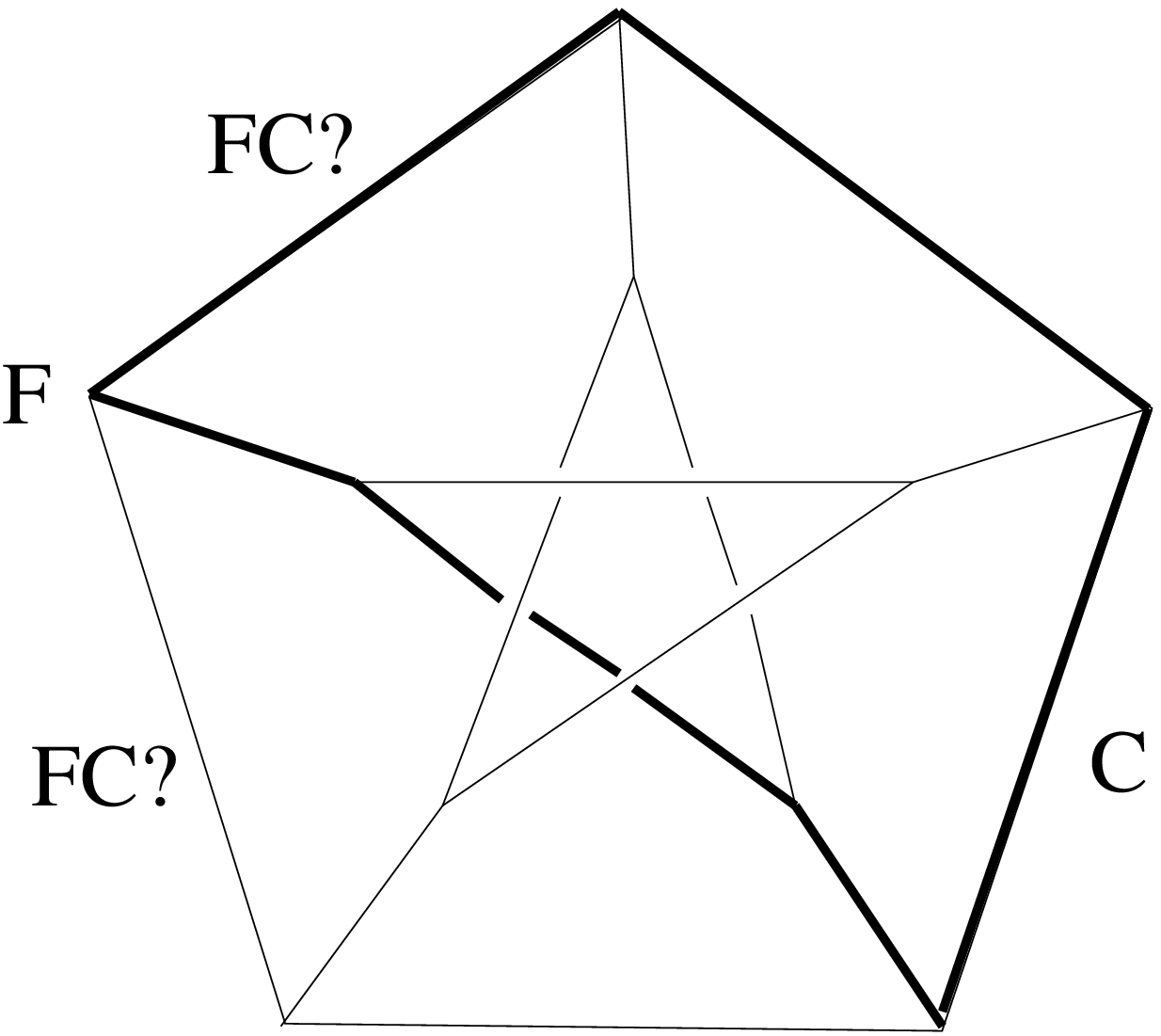,height=3cm,width = 3cm}}
\end{tabular}
\caption{The Petersen graph - an almost-building?}
\label{petersen}
\end{figure}

\noindent
The Petersen graph $\D$:
Consider a node $v$ with
5 labelled leaves 
\epsfig{file=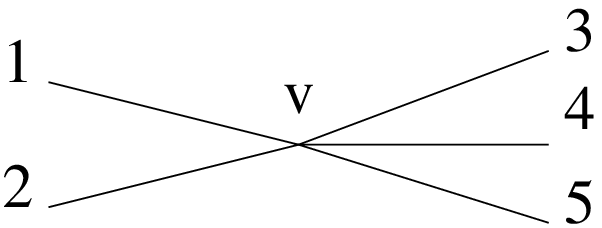,height=0.5cm,width = 1.2cm}.
This corresponds to the empty face of $\D$. 
To get the rest of $\D$, we do a ``blowup'' at $v$.
The trees of the form \epsfig{file=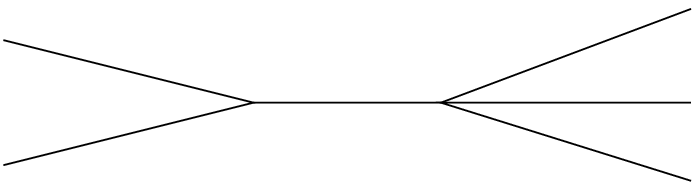,height=0.4cm,width = 1.2cm}
(resp. \epsfig{file=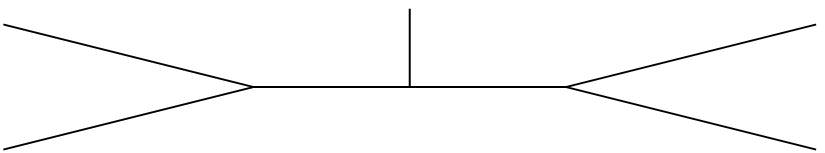,height=0.4cm,width = 1.3cm} )
are the vertices (resp. edges) of $\D$.
The 5 leaves of each tree are labelled $1,2,3,4$ and $5$ in some order.
Define an apartment to be the subcomplex made of those edges (chambers)
that have the form \epsfig{file=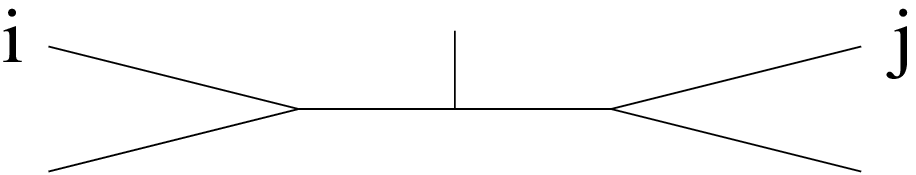,height=0.5cm,width = 1.2cm}
for a fixed $i,j$.
There are 10 apartments, 15 chambers and every chamber lies in exactly 4
apartments. 
Furthermore each apartment is an hexagon and 
given two chambers, there is always an apartment $\Sigma$ containing
them. So we may first define projection maps in $\Sigma$ as usual and then
try to extend them to $\D$.
However, in some cases, there are two apartments containing a 
face $F$ and a chamber $C$ and the products are not always consistent.
This happens exactly when $F$ and $C$ are opposite to each other in
some pentagon as shown in Figure~\ref{petersen}. Then there are exactly two
apartments containing $F$ and $C$ one of which is indicated in the figure.
And the two candidates for $FC$ are indeed the two chambers containing
$F$ that are closest to $C$.

$(ii)$
In~\cite{\bjor}, Bj{\"o}rner
constructs many examples of shellable complexes $\D$.
For example, he shows that all finite semimodular and supersolvable
lattices are shellable.
As a very special example, we would also like to mention the Whitehead
poset which has found an application to group theory \cite{\brady}.
In all these examples, it is clear that $\D$ has many shellings.
It would be worthwhile to understand the precise compatibility
relations among these shellings.

$(iii)$
New examples from old:
Let $\D$ be a labelled simplicial complex that satisfies our axioms.
Then what can be said about its type
selected subcomplex $\D_J$ or its $k$-skeleta $\D_k$ (if $\D$ is not
labelled)? We may ask the same question for the barycentric
subdivision $\D^{\prime}$. These complexes again have many shellings
but they do not directly give us new examples.
However they might point to an appropriate generalisation of the
theory.

$(iv)$ Finally, we mention non-simplicial examples like polytopes and
hyperplane arrangements.  The second situation may be easier to handle
since we already have projection maps in that case.
The other non-simplicial examples to consider are LRBs
as explained in Section~\ref{s:LRB}.
We mention that there is a notion of
shelling for a poset. It 
is called a recursive coatom ordering~\cite{\red}. 
This would be the relevant shelling concept in this generality.

\vanish{
\begin{figure}[hbt]
\centering
\begin{tabular}{c@{\qquad}c}
\mbox{\epsfig{file=round.eps,height=3cm,width = 3cm}}
\end{tabular}
\caption{} 
\label{round}
\end{figure}

\begin{example}
If $\D$ is an even-sided polygon, then we can modify
Example~\ref{e:nmetric} slightly so that $\D$ satisfies the opposite
axiom. The situation is easiest to describe from the point of view of
shellings. For any $C \in \C$, the partial order $\leq_C$ is a total
order beginning at $C$ and going in either the clockwise or the
anti-clockwise direction. In Example~\ref{e:nmetric}, they all went in
the same direction. Now we require that for adjacent chambers the
directions are opposite. See Figure. 
\end{example}
}
\subsection*{Acknowledgments} 
I would like to thank my advisor Ken Brown
for his encouragement, sharing of his insights
and patient reading of earlier versions of this paper.
I would also like to thank 
H. Abels, P. Abramenko, M. Aguiar, L. Billera,
A. Bj{\"o}rner, R. Forman, J. Meier, B. M{\"u}hlherr and
G. Ziegler for useful discussions, comments and correspondence.



\providecommand{\bysame}{\leavevmode\hbox to3em{\hrulefill}\thinspace}

\end{document}